\documentclass[12pt]{article}

\usepackage[cp1251]{inputenc}
\usepackage[english]{babel}

\usepackage{amsfonts,epsfig,epstopdf}
\usepackage{amssymb,amsmath,enumerate,float}
\makeatother

\newcommand{\labell}[1]{%
\label{#1}}

\addto\captionsenglish{}

\renewcommand{\i}{i}
\newcommand{\ptl}{\partial}

\newcommand{\be}{\begin{equation}}
\newcommand{\ee}{\end{equation}}
\newcommand{\beq}{\begin{equation}}
\newcommand{\eeq}{\end{equation}}

\newcommand{\bs}{\mathcal{R}}
\newcommand{\eps}{\varepsilon}

\title{Diffraction of a mode close to its cut-off by a transversal screen in a planar waveguide}

\author{A. V. Shanin, A. I. Korolkov}

\begin{document}

\maketitle

\begin{abstract}

The problem of diffraction of a waveguide mode by a thin Neumann screen
is considered. The incident mode is assumed to have frequency close to the
cut-off. The problem is reduced to a propagation problem on a branched
surface and then is considered in the parabolic approximation. Using the
embedding formula approach, the reflection and transmission coefficients
are expressed through the directivities of the edge Green's function of the
propagation problem. The asymptotics of the directivities of the edge Green's
functions are constructed for the case of small gaps between the screen and the
walls of the waveguide. As the result, the reflection and transmission coefficients are
found. The validity of known asymptotics of these coefficients is studied.

\end{abstract}

\section{Problem formulation and introductory notes}

\begin{figure}[ht]
\centerline{\epsfig{file=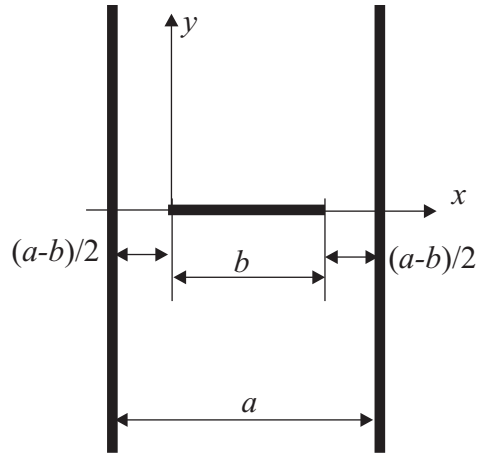}}
\caption{Geometry of the waveguide}
\label{fig01}
\end{figure}

Consider a planar acoustic waveguide in the plane $(x,y)$ composed of two acoustically hard walls
located at $x= (b-a)/2$ and $x = (b+a)/2$ (see Fig.~\ref{fig01}). The width of the
waveguide is $a$; the position of the walls is chosen for convenience of computations.
The Helmholtz equation
\begin{equation}
\Delta \tilde u + k^2 \tilde u = 0
\labell{eq0101}
\end{equation}
is fulfilled inside the waveguide by the total field $\tilde u$.
The time dependence of the form $e^{-i\omega t}$ is omitted
everywhere. The wavenumber $k$ has a small positive imaginary part
corresponding to absorption in the medium. This feature enables us
to avoid discussing trapped modes and to consider all series in the paper as convergent.
A thin acoustically hard (Neumann) screen is located inside the waveguide, namely
it occupies the segment $y = 0$, $0< x < b$. The width of the
screen is equal to $b$. Of course $b < a$.

Throughout the paper we use notations with tilde for values related to the Helmholtz equation and
notations without tilde for values related to the parabolic equation.

An incident  waveguide mode falls from the domain $y > 0$ and is scattered by the screen.
The mode has form
\begin{equation}
\tilde u_{\rm in} = \cos \left( k (x+(a-b)/2) \cos \tilde \theta_{\rm in} \right)
e^{ -i k y \sin \tilde \theta_{\rm in} }
\labell{eq0102}
\end{equation}
for some $\tilde \theta_{\rm in}$ obeying the relation
\begin{equation}
\theta_{\rm in} = \tilde \theta_m.
\labell{eq0103}
\end{equation}
The angles $\theta_n$, $n \in \mathbb{Z}$, correspond to waveguide modes:
\begin{equation}
ka \cos \tilde \theta_n = \pi n.
\labell{eq0104}
\end{equation}
Parameter $m$ (the index of the incident mode) is assumed to be fixed.

The scattered field $\tilde u_{\rm sc}$
can be represented as a linear combination of waveguide modes. For
$y> 0 $ the scattered field is
$\tilde u_{\rm sc} \equiv \tilde u - \tilde u_{\rm in}$,
and
the expansion is as follows:
\begin{equation}
\tilde u_{\rm sc} (x,y) = \sum_{n=0}^{\infty} \tilde R_{n,m}
\cos \left( k (x+(a-b)/2) \cos \tilde \theta_n \right)
e^{ i k y \sin \tilde \theta_n }.
\labell{eq0105}
\end{equation}
For $y<  0 $ we define $\tilde u_{\rm sc} \equiv \tilde u $, and the expansion is as follows:
\begin{equation}
\tilde u_{\rm sc} (x,y) = \sum_{n=0}^{\infty} \tilde T_{n,m}
\cos \left( k (x+(a-b)/2) \cos \tilde \theta_n \right)
e^{ - i k y \sin \tilde \theta_n },
\labell{eq0105a}
\end{equation}
Coefficients $\tilde R_{n,m}$ and $\tilde T_{n,m}$ are elements of the reflection and transmission matrices.
Note that $R_{n,m} = 0$ and $T_{n,m} = 0$ if $n-m$ is odd because of the symmetry of the problem.
Our aim is to find $\tilde R_{n,m}$ and $\tilde T_{n,m}$.

By applying the
principle of reflections one can reduce the problem of scattering in the waveguide to the problem of
diffraction by an infinite periodic grating composed of Neumann segments (see below). It is well-known that
such a problem leads to a $2\times 2$ matrix Wiener--Hopf problem. The Wiener--Hopf problem remains unsolved
except the only particular case of $b = a/2$ \cite{Lukianov1980,Lukianov1981, Daniele1990, Abrahams2007}, for which it reduces to a matrix problem with matrices of the Daniele--Khrapkov form. The latter can be factorized explicitly. For details see survey articles \cite{Luneburg1993, Luneburg2004}. For normal incidence this problem has been solved in \cite{Baldwin1954,Weinstein1969} by scalar Wiener--Hopf technique using dual integral equation formulation.
For the case of arbitrary incidence this problem has been solved first by Riemann--Hilbert method~\cite{Luneburg1971}.

In \cite{Achenbach1986,Sumbatyan1996}  the problem is reduced to a singular integral equation, which is solved approximately. In \cite{Porter1996} it is considered with help of eigenfunction expansion technique.

Another method that can be applied to the problem is the method of matching series expansions
\cite{Agranovich1962, Shestopalov1971}. For this method the values $\tilde R_{n,m}$, $\tilde T_{n,m}$ are unknowns of an infinite
(truncated somehow) linear algebraic system.

In the current paper we are not presenting any method for solving a general problem of diffraction by
a transversal screen in a waveguide. Instead, we study asymptotically the case of
\begin{equation}
\tilde \theta_{\rm in} \ll 1,
\labell{eq0106}
\end{equation}
i.\ e.\ we assume that
partial Brillouin waves propagate almost normally to the axis of the waveguide.
Physically, in this case the time frequency of the incident wave is close to the cut-off
frequency of the $m$-th mode.
Condition (\ref{eq0106}) leads to some
interesting properties of the scattering process. A detailed
theory can be found in \cite{Nazarov2015,Nazarov2016}. Depending on the structure of the scatterer,
there can occur either a general case of an almost complete
reflection ($\tilde R_{m,m} \to -1$ as $\tilde \theta_{\rm in} \to 0$),
or a particular case of an anomalous transmission.  According to the classification given in \cite{Nazarov2016}, in the case of a thin
transversal hard screen one should observe the anomalous transmission with
$\tilde T_{m,m} \to 1$ as $\tilde \theta_{\rm in} \to 0$.
The classification of \cite{Nazarov2016} is based on the existence of the so-called {\em threshold stabilizing solutions}. Namely, it is proven that an anomalous transmission takes place if and only if there exists
a solution of the scattering problem having
a special form at the cut-off (threshold) frequency of the $m$-th mode. This solution should
have structure of the non-growing $m$-th waveguide mode (plus, possibly, an exponentially decaying remainder) in the branches of the waveguide. In our case such a solution can be easily written:
\[
\tilde u^{\rm st} = \cos \left( k (x+(a-b)/2) \cos \tilde \theta_m \right).
\]
The solution is constant with respect to $y$, so it obeys Neumann boundary conditions on the hard screen.

The anomalous transmission is a complicated and surprising phenomenon.
It is a result of a resonance interaction of gaps between the screen and each of the walls.
The phenomenon is formulated as a mathematical theorem, and the theorem is valid for any $b$
smaller than~$a$.
According to the theorem,
it should be $\tilde T_{m,m} \to 1$
for small but fixed $(a-b)/2$ and $\tilde \theta_{\rm in} \to 0$.
At the same time, for small but fixed $\tilde \theta_{\rm in}$ and $(a-b)/2 \to 0$ it should be $\tilde T_{m,m} \to 0$
(since the gap becomes small and the wave no longer can  penetrate it).
Therefore, in the domain of two small parameters $\tilde \theta_{\rm in}$ and $(a-b)/a$ there should be a
boundary separating these two regimes. The aim of the current paper is to study this boundary and to describe the regimes in more details. The problem has been formulated by S.A.Nazarov, and the authors are grateful to him for
inspiring discussions.

We should note that a total (resonance) transmission
through a periodic grating has been described by Maljuzhinets for a different physical situation.
The work of Maljuzhinets was dedicated to creation of acoustical environment of Palace of Soviets in Moscow
(was not constructed). The phenomenon of total transmission is called now the {\em  Maljuzhinets phenomenon}.

The structure of the paper is as follows.
In Section~\ref{sec02} the problem of scattering by a thin hard wall in a waveguide is transformed into the
problem of propagation on a branched surface $\bs$. The principle of reflection is used for this. The
resulting branched surface has two sheets and an infinite number of branch points of order~2. Formulae linking the
transmission and reflection coefficients for the waveguide problem and those related to the branched surface problem
are derived.

In Section~\ref{sec03} the Helmholtz equation on surface $\bs$ is reduced to a parabolic equation. The assumptions of  shortness of the wave and of smallness of the angle of incidence are used for this.
Again,  transmission and reflection coefficients are introduced for the parabolic problem. These coefficients
are linked with the waveguide reflection and transmission coefficients by simple relations.

In Section~\ref{sec04} edge Green's functions for the parabolic problem on the branched surface are introduced. The
edge Green's functions are generated by point sources located near the branch points rather than by an incident plane wave.
In Section~\ref{sec05} the Green's theorem for the parabolic equation is derived.
Using this theorem the directivities
of the edge Green's functions are expressed through integrals of the wave fields.

In Section~\ref{sec06} the edge values of the diffraction field are introduced. They are the values of the wave field at the
branch points. The edge values play an important role in the consideration. A formula connecting the edge values
of the wave field with the directivities of the edge Green's function is derived. This formula is based on the
reciprocity principle, which is proven to be valid for the parabolic equation on branched surfaces.

In Section~\ref{sec07} an embedding formula for the field is derived. The embedding formula is an expression
connecting the diffracted wave field (generated by the incident plane wave) with the edge Green's functions.
A standard technique based on an embedding operator is used for derivation. Unfortunately, the embedding formula for the
field is not enough for our consideration. We need expressions for the transmission and reflection coefficients,
i.\ e.\ the embedding formula for the coefficients. Sections~\ref{sec08}, ~\ref{sec09}, and~\ref{sec10}
are dedicated to this.
Formulae in Sections~\ref{sec08} and~\ref{sec09} are obtained in a standard way. However, these formulae do not enable one to
obtain the zero-order transmission coefficient since the embedding operator nullifies corresponding wave. Thus,
the formula from Section~\ref{sec09} should be regularized. A regularization is built in Section~\ref{sec10}. This regularization is
based on the concept of extended directivity of an edge Green's function. This concept requires an alternative formulation of the diffraction problem developed in Appendix~A.

Sections~\ref{sec11} and~\ref{sec12} are dedicated to estimation of the directivities of the edge Green's functions. The consideration is based on the model of scattering by a short Dirichlet or Neumann segment. Such a model is derived in Appendix~B.

In Section~\ref{sec13} some analytical and numerical results are presented. The boundary between transmission
and reflection regimes is found.
This boundary is given by formulae (\ref{eq1301}) and (\ref{eq1302}).


\section{Preparatory steps}
\subsection{Application of the reflection principle. A problem on a branched surface}
\label{sec02}

Apply the reflection principle to the walls of the waveguide, i.\ e.\ consider the field $\tilde u(x,y)$
on the whole plane $(x,y)$ instead of the strip
\[
-(a-b)/2< x < (a+b)/2 .
\]
Since the walls are hard, we use an
even continuation of the field through the walls:
\begin{equation}
\tilde u (x, y) = \tilde u (b-a-x , y),
\qquad
\tilde u (x, y) = \tilde u (b+a-x , y).
\labell{eq0201}
\end{equation}
Now one can study diffraction of wave (\ref{eq0102}) by a periodic set of infinitely thin acoustically
hard screens (segments $y = 0$, $an< x < an + b $, $n \in \mathbb{Z}$).

\begin{figure}[ht]
\centerline{\epsfig{file=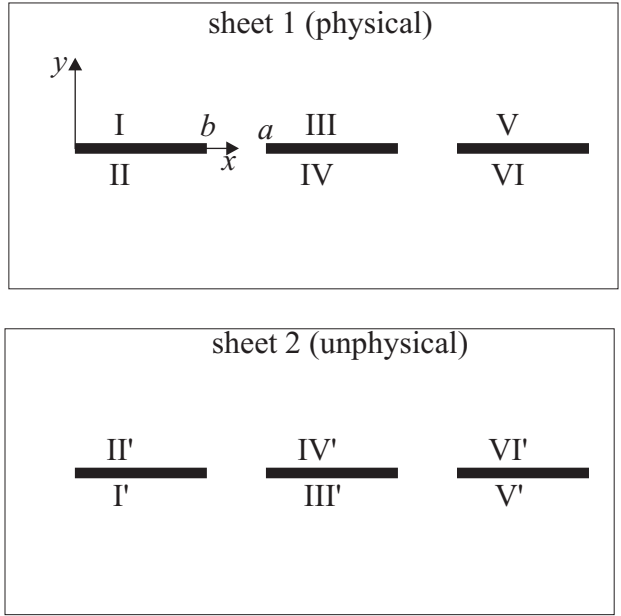}}
\caption{Scheme of the branched surface}
\label{fig02}
\end{figure}

Apply the reflection principle to the hard screens. Namely, cut the plane along the screens and duplicate the
plane with the cuts. On the second (unphysical) sheet define an even reflection of the wave field
by the formula
\begin{equation}
\tilde u'(x,y) = \tilde u(x,-y).
\labell{eq0202}
\end{equation}
Attach the shores of the cuts of the unphysical sheet to those of the physical sheet as it is
shown in Fig.~\ref{fig02}. Shores labeled by the same Roman numbers should be linked together. As the result, get
a branched surface $\bs$ having two sheets and an infinite number of branch points. Each branch point has order~2.
The positions of the branch points are $x_n$:
\begin{equation}
x_{2n} = a n, \qquad x_{2n + 1} = an + b , \qquad n \in \mathbb{Z}.
\labell{eq0203}
\end{equation}
The field should obey  Meixner's conditions in the vicinities of the branch points
(the energy of the field should be locally finite there).

Now instead of scattering by hard screens in a plane we will study propagation on $\bs$.
The incident field is comprised of two plane waves on~$\bs$.
 One is (\ref{eq0102}) existing in the upper half-plane
on the physical sheet. The other is $\tilde u_{\rm in}(x,-y)$ existing in the lower half-plane of the unphysical
sheet.

\begin{figure}[ht]
\centerline{\epsfig{file=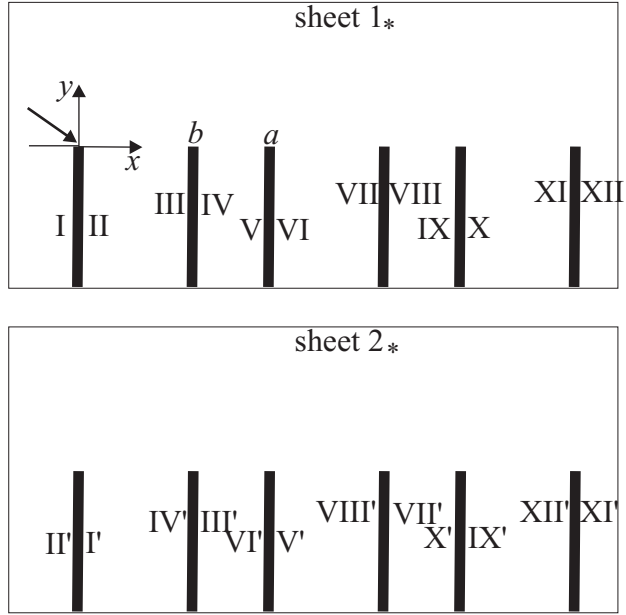}}
\caption{An alternative representation of $\bs$}
\label{fig03}
\end{figure}

There is an alternative representation of $\bs$. Namely, the cuts can be drawn along the
half-lines $x = x_n$, $y < 0$. The scheme of the surface is shown in Fig.~\ref{fig03}.
The upper half-plane of  sheet~$1_*$ coincides with the upper half-plane of sheet~1, the same is valid for the
upper half-planes of sheets 2  and~$2_*$. However, the lower half-planes are cut differently in these representations.

To avoid misunderstanding we will denote the field on the sheets cut as in Fig.~\ref{fig02}
by $\tilde u(x,y,\nu)$, where $\nu = 1,2$ is the index of the sheet. If the surface is cut according to Fig.~\ref{fig03} then the field will be denoted as
$\tilde u(x,y,\nu_*)$, $\nu = 1,2$. Formally, the star belongs not to $\nu$, but to the notation
$\tilde u(\cdot,\cdot, \cdot_*)$

The cuts shown in Fig.~\ref{fig02} and Fig.~\ref{fig03} are
linked by the relation
\begin{equation}
\tilde u(x,y,\nu) = \tilde u(x,y,\nu_*) \qquad \mbox{for} \quad y>0
\quad \mbox{or} \quad y < 0, \, \, x_{2m+1}<x<x_{2m+2}
\labell{eq0203a}
\end{equation}
\begin{equation}
\tilde u(x,y,\nu) =
\tilde u(x,y,(3-\nu)_*) \qquad \mbox{for} \quad y<0,\, \, x_{2m}<x<x_{2m+1}
\labell{eq0203b}
\end{equation}

The problem of propagation on $\bs$ is slightly more general than the problem of scattering
by a hard wall in a waveguide. Namely, instead of two incident cosine  waves $\tilde u_{\rm in}(x,y)$,
$\tilde u_{\rm in}(x,-y)$
existing on two sheets one can consider a single incident exponential wave
\begin{equation}
\bar u_{\rm in} = e^{ i k (x \cos \tilde \theta_{\rm in} -y \cos \tilde \theta_{\rm in} )}
\labell{eq0204}
\end{equation}
which is present on the upper half-plane ($y > 0$) of sheet~1.
Moreover, it is important that for the problem on $\bs$ the parameter $\tilde \theta_{\rm in}$
is an arbitrary angle, not necessary obeying (\ref{eq0103}), (\ref{eq0104}).

Consider the problem formulated above on $\bs$. Let the incident wave (\ref{eq0204}) be falling
along the first sheet. Let $\tilde \theta_{\rm in}$ {\em not obey} (\ref{eq0103}), (\ref{eq0104}).
Denote the scattered field by
$\bar u_{\rm sc}(x,y,\nu)$.
The scattered field is connected with the total field $\bar u$ by the relation
\begin{equation}
\bar u(x,y,1) = \bar u_{\rm sc}(x,y,1) + \bar u_{\rm in}(x,y)\quad \mbox{ for }y>0,
\labell{eq0204a}
\end{equation}
\[
\bar u(x,y,1) = \bar u_{\rm sc}(x,y,1) \quad \mbox{ for }y < 0,
\]
\[
\bar u(x,y,2) = \bar u_{\rm sc}(x,y,2).
\]

According to Floquet principle, the scattered field can be represented as series. For $y > 0$
the series take form
\begin{equation}
\bar u_{\rm sc}(x,y, \nu) = \sum_{n = -\infty}^{\infty} \bar R_n^\nu
e^{ i k (x \cos \tilde \psi_n +y \cos \tilde \psi_n )}
\labell{eq0205}
\end{equation}
for some coefficients $\bar R_n^\nu$, while for $y < 0$ the series take form
\begin{equation}
\bar u_{\rm sc}(x,y, \nu) = \sum_{n = -\infty}^{\infty} \bar T_n^\nu
e^{ i k (x \cos \tilde \psi_n -y \cos \tilde \psi_n )} .
\labell{eq0206}
\end{equation}
Here
\begin{equation}
\tilde \psi_n = \arccos ( \cos \tilde \theta_{\rm in} - 2 \pi n/(ka)).
\labell{eq0207}
\end{equation}

The ``old'' field $\tilde u (x,y)$ can be represented through the ``new'' values
$\bar u(x,y,\nu)$ as follows. Consider the limit
\begin{equation}
\tilde \theta_{\rm in} \to \tilde \theta_m.
\labell{eq0208}
\end{equation}
Obviously,
\begin{equation}
\tilde u(x,y) = \frac{1}{2} e^{ ik\cos \tilde \theta_{\rm in} (a-b)/2 } \times
\labell{eq0209}
\end{equation}
\[
\left( \bar u(x, y, 1) + \bar u(x, -y, 2) + (-1)^m \bar u(b-x, y, 1) + (-1)^m \bar u(b-x, -y, 2)
\right)
\]
and
\begin{equation}
\tilde T_{m-2n,m} = e^{ ik( \cos \tilde \theta_{\rm in} - \cos \tilde \theta_{m+2n}) (a-b)/2 }
\times
\labell{eq0210}
\end{equation}
\[
(\bar T_n^1 + \bar R_n^2 + \bar T_{-n-m}^1 + \bar R_{-n-m}^2),
\]
\begin{equation}
\tilde R_{m-2n,m} = e^{ ik( \cos \tilde \theta_{\rm in} - \cos \tilde \theta_{m+2n}) (a-b)/2 }
\times
\labell{eq0211}
\end{equation}
\[
(\bar R_n^1 + \bar T_n^2 + \bar R_{-n-m}^1 + \bar T_{-n-m}^2).
\]

Write down one more elementary properties of the coefficients $\bar T$, $\bar R$.
Consider the field $\bar u(x,y,1) + \bar u(x,y,2)$. This field obeys the Helmholtz equation on the
$(x,y)$--plane without branch points. Thus, it is the trivial solution comprised of the incident plane wave.
Therefore
\begin{equation}
\bar R^1_n + \bar R^2_n  = 0 ,
\qquad
\bar T^1_n + \bar T^2_n  = \delta_{n, 0} ,
\labell{eq0212}
\end{equation}
where $\delta$ is the Kronecker's delta.


\subsection{Parabolic equation}
\label{sec03}

Consider the problem for $\bar u$ on $\bs$.
To simplify the consideration study the short-wave case, i.\ e.\ the values $a$, $b$ and $(a-b)/2$ are large
comparatively to the wavelength.
Moreover, assume that $\tilde \theta_{\rm in} \ll 1$. In this case
one can apply the parabolic approximation of diffraction theory \cite{Talanov1995,Fock,Weinstein}
taking the positive $x$-direction as the main propagation direction.
Within this approximation the
field on $\bs$ is represented in the form
\begin{equation}
\bar u(x,y) = e^{i k x} u(x,u).
\labell{eq0301}
\end{equation}
Dependence of $u(x,y)$ on $x$ and $y$ is assumed to be slower than that of the exponential factor on~$x$.
The Helmholtz equation (\ref{eq0101}) is approximated by the parabolic equation
\begin{equation}
L[u] = 0,
\qquad
L \equiv
\ptl_x + (2 i k)^{-1} \ptl_y^2.
\labell{eq0302}
\end{equation}
The parabolic equation describes well the Fresnel diffraction (when the angle of scattering is small),
describes badly diffraction at large angles, and it neglects
the back-scattering (there are no waves traveling in the negative $x$-direction in this approximation).
Fortunately, the problem of interest, i.~e.\ the diffraction of a waveguide wave close to the cut-off
is mainly described as the Fresnel scattering.

It is known that the parabolic approximation works well for the angles $\tilde \theta$ for which the
dimensionless combination
$\sqrt{ka} \, \tilde \theta$ is not large compared to~1. Thus, we demand that
$\sqrt{ka} \, \tilde \theta_{\rm in}$ is not large. Also we expect that the coefficients
$\bar T_{n,m}$, $\bar R_{n,m}$ will be estimated reasonably only for $n$
such that $\sqrt{ka} \, \tilde \psi_n$ is not large.

The incident wave (existing only in the upper half-plane of sheet~1) has form
\begin{equation}
u_{\rm in} = e^{ -ikx \theta_{\rm in}^2/2 - i k \theta_{\rm in} y },
\labell{eq0303}
\end{equation}
where $\theta_{\rm in}$ is a parameter approximately equal to $\tilde \theta_{\rm in}$, such that
\begin{equation}
1-\frac{\theta_{\rm in}^2}{2} = \cos \tilde \theta_{\rm in}.
\labell{eq0304}
\end{equation}
One can see that in some strip $|y|  < \mbox{const}$ the wave (\ref{eq0303}) together
with (\ref{eq0301}) provides an approximation to $\tilde u_{\rm in}$.

We look for a scattered field $u_{\rm sc}(x,y,\nu)$ on $\bs$. The scattered field is linked with
the total field $u$ in a way similar to (\ref{eq0204a}):
\begin{equation}
u(x,y,1) = u_{\rm sc}(x,y,1) + u_{\rm in}(x,y)\quad \mbox{ for }y>0,
\labell{eq0304a}
\end{equation}
\[
u(x,y,1) = u_{\rm sc}(x,y,1) \mbox{ for }y < 0,
\]
\[
u(x,y,2) = u_{\rm sc}(x,y,2).
\]
The total field should be
continuous on $\bs$ (except the branch points) and it should be bounded near the branch points.
The last condition plays the role of Meixner's condition and guarantees
the absence of sources at the branch points.

According to the Floquet principle, the scattered field should be represented as series
\begin{equation}
u_{\rm sc}(x,y,\nu) = \sum_{n=-\infty}^{\infty} R^\nu_n
e^{ -ikx \psi_n^2/2 + i k  \psi_n y }
\labell{eq0305a}
\end{equation}
for $ y > 0$, and
\begin{equation}
u_{\rm sc}(x,y,\nu) = \sum_{n=-\infty}^{\infty} T^\nu_n
e^{ -ikx \psi_n^2/2 - i k  \psi_n y}
\labell{eq0305b}
\end{equation}
for $y < 0$. These representations are parabolic approximations of (\ref{eq0205}) and
(\ref{eq0206}). The ``angles'' $\psi_n$ are defined by
\begin{equation}
\psi_n = \sqrt{\theta_{\rm in}^2 + \frac{4 \pi n}{ka}}.
\labell{eq0306}
\end{equation}
Note that $\psi_0 = \theta_{\rm in}$.

Obviously, for not large $|n|$,  $\psi_n \approx \tilde \psi_n$, and
\begin{equation}
\bar T^\nu_n \approx T^{\nu}_n,
\qquad
\bar R^\nu_n \approx R^{\nu}_n.
\labell{eq0307}
\end{equation}
Moreover, we neglect the backscattering, which corresponds to the last two terms in the sums
in (\ref{eq0210}) and (\ref{eq0211}).
Thus, we get from (\ref{eq0210}), (\ref{eq0211}), and (\ref{eq0307})
\begin{equation}
\tilde T_{m-2n,m} \approx e^{ i \pi n (a-b)/a} ( T_n^1 +  R_n^2 ),
\labell{eq0308}
\end{equation}
\begin{equation}
\tilde R_{m-2n,m} \approx e^{ i \pi n (a-b)/a} ( R_n^1 +  T_n^2 ).
\labell{eq0309}
\end{equation}
The last formulae connect the solution of the parabolic problem with that of the initial Helmholtz
problem in the waveguide.

Similarly to (\ref{eq0212}),
\begin{equation}
 R^1_n +  R^2_n  = 0 ,
\qquad
 T^1_n +  T^2_n  = \delta_{n, 0} .
\labell{eq0310}
\end{equation}

The parabolic equation of diffraction theory has one important property making the
parabolic problems much simpler than their Helmholtz prototypes. Namely, if the strip
$x' < x < x''$ is free of branch points, then for any solution $v$ of equation (\ref{eq0302})
in this strip
\begin{equation}
v(x'', y) = \int \limits_{-\infty}^{\infty} v(x',y') G(x''-x' , y-y') dy',
\labell{eq0311}
\end{equation}
where
\begin{equation}
G(x,y) = \sqrt{\frac{k}{2 \pi i x}} e^{ i k y^2/(2 x)}
\labell{eq0312}
\end{equation}
is the Green's function of the whole plane, i.\ e.\  a solution of the equation
\[
L[G] = \delta(x) \delta(y)
\]
equal to zero for $x < 0$.

Let $v$ be defined on $\bs$, decay as $|y|\to \infty$, and be bounded near the branch points.
Let be $x_n<x<x_{n+1}$. The integral representation (\ref{eq0311})
can be rewritten as follows:
\begin{equation}
v(x, y, \nu_*) =  \int \limits_{0}^{\infty} v(x_n,y', \nu_*) G(x-x_n , y-y') dy' +
\int \limits_{-\infty}^0 v(x_n+0,y', \nu_*) G(x-x_n , y-y') dy' =
\labell{eq0313}
\end{equation}
\[
\int \limits_{0}^{\infty} v(x_n,y', \nu_*) G(x-x_n , y-y') dy' +
\int \limits_{-\infty}^0 v(x_n-0,y', (3-\nu)_*) G(x-x_n , y-y') dy'
\]
Note the usage of $x_n+0$ and $x_n-0$ notations in (\ref{eq0313}). Since the surface is cut into sheets~$1_*$
and~$2_*$ along the lines $x= x_n$, these notations indicate whether the values on the right or on the
left shore of a cut are taken.

Below we use the principle of limiting absorption to separate the outgoing and the ongoing
waves. As we noted before, absorption corresponds to $k$ having a small positive
imaginary part,
i.\ e.\
\[
k = k' (1+i\epsilon), \qquad \epsilon \to 0, \qquad {\rm Im}[k'] =0.
\]
At the same time, we keep the Floquet constant
$
e^{ -i k a \theta_{\rm in}^2 /2 }
$
having absolute value equal to~1, thus $k \theta_{\rm in}^2$ should be real. Therefore
\[
\theta_{\rm in} = \theta'_{\rm in} (1 - i \epsilon /2), \qquad {\rm Im}[\theta'_{\rm in}] =0 .
\]
The combination $k \theta_{\rm in}$ has a {\em positive} imaginary part:
\[
k \theta_{\rm in} \approx k' \theta'_{\rm in} (1 + i \epsilon / 2).
\]
A similar reasoning shows that for each $\psi_n$ the combination $k \psi_n$, $n \in \mathbb{Z}$
has a positive imaginary part (vanishing or not), while $k \psi_n^2$ is real. Thus, the ``plane wave''
\[
e^{ -ikx \psi_n^2/2 + i k \psi_n y  }
\]
for any $n$ oscillates in $x$-direction and decays in positive $y$-direction.


\section{Edge Green's functions}
\subsection{Introduction of the edge Green's functions}
\label{sec04}

Edge Green's functions are the fields generated by sources located near the branch points.
These functions are used for derivation of embedding formulae \cite{Shanin}.

\begin{figure}[ht]
\centerline{\epsfig{file=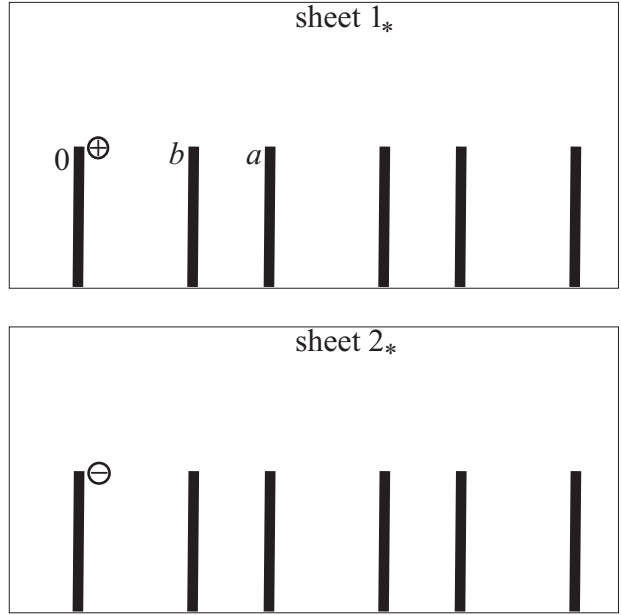}}
\caption{Sources for the edge Green's function $w_0$}
\label{fig04}
\end{figure}

Introduce edge Green's functions $w_n (x,y,\nu_*)$, where $n \in \mathbb{Z}$ is the number of the branch point,
the notation of the sheets is as above.
Each edge Green's function is generated by two sources located to the right of the point $x_n$:
\begin{equation}
L[w_n](x,y,1_*) = \delta(x-(x_n + 0)) \delta(y) ,
\labell{eq0401}
\end{equation}
\begin{equation}
L[w_n](x,y,2_*) = -\delta(x-(x_n + 0)) \delta(y) ,
\labell{eq0402}
\end{equation}
where $\delta$ is the Dirac delta-function.
The location of the sources for the edge Green's function $w_0$ is shown in Fig.~\ref{fig04}.

The notation $x_n+0$ indicates that the source is located to the right of $x_n$. According to general
properties of the parabolic equation, in this case
$w(x,y,\nu_*)= 0$ for $x < x_n$, and
\[
w_n (x,y,\nu_*) = (-1)^{\nu-1} G(x-x_n , y)
\]
for $x_n<x<x_{n+1}$. Since the source is located to the right of $x_n$, the field does not feel the
cut at $x = x_n$.

The structure of the sources and the symmetry of surface $\bs$ lead to the relations
\begin{equation}
w_{n}(x,y,1_*) = - w_n(x,y,2_*),
\labell{eq0403a}
\end{equation}
\begin{equation}
w_{n}(x,-y,1) = (-1)^{n+1} w_n(x,y,1).
\labell{eq0403}
\end{equation}

Consider edge Green's function $w_1$ in the upper half-plane of sheet~$1$. Similarly to the
Helmholtz case, far from the source the field can be asymptotically
represented as a Green's function of the plane
multiplied by a factor depending only on the angle of scattering:
\begin{equation}
w_n(x,y,1_*) = V_n \left( \frac{y}{x-x_n} \right) G(x-x_n,y) + o(G(x-x_n,y)).
\labell{eq0501}
\end{equation}
This representation defines uniquely directivity $V_n(\theta)$.
Directivity $V_n(\theta)$ is defined on the half-line $\theta > 0$.
Since $\theta = 0$ corresponds to the $x$-axis, we don't expect that $V_n(\theta)$ can be smoothly
continued through $\theta = 0$.

Definition (\ref{eq0501}) can be used in the upper half-plane of sheet~$1$, or, the same, of sheet~$1_*$.
Using the notation of sheets of the surface, we can write $V_n(\theta) = V_n(\theta, 1)$, $\theta > 0$.
The directivities at the other three half-planes can be found from the symmetry relation (\ref{eq0403}):
\begin{equation}
V_n(\theta,2) = - V_n(\theta), \qquad \theta > 0,
\labell{eq0502}
\end{equation}
\begin{equation}
V_n(\theta,\nu) = (-1)^{\nu+n} V_n(-\theta), \qquad \theta < 0.
\labell{eq0503}
\end{equation}

Since $\bs$ is periodic in the $x$-direction, for $l \in \mathbb{Z}$
\begin{equation}
w_{n+2l} (x+a l , y,\nu) = w_n(x,y,\nu)
\labell{eq0504}
\end{equation}
\begin{equation}
V_{n+2l} (\theta)= V_n(\theta).
\labell{eq0505}
\end{equation}
Thus, it is sufficient to find only $V_0$ and $V_1$ to know all directivities of the edge Green's functions.


\subsection{Green's theorem and an integral representation for the directivities of the edge Green's
functions}
\label{sec05}

Derive Green's theorem for the parabolic equation.
Introduce  operator
\begin{equation}
L' = - \ptl_x + (2ik)^{-1} \ptl^2_y.
\labell{eqA101}
\end{equation}
Let $v$ and $w$ be some functions obeying equations
\begin{equation}
L[v] = f(x,y),
\qquad
L'[w] = h(x,y).
\labell{eqA102}
\end{equation}
in some domain $\Omega$ of a plane having piecewise--linear boundary $\ptl \Omega$.
Then
\begin{equation}
\frac{1}{2 i k}\int \limits_{\ptl \Omega} [({\bf v}\cdot {\bf n})w - ({\bf w}\cdot {\bf n})v] dl =
\int \limits_\Omega [f w - h v] ds.
\labell{eqA103}
\end{equation}
where ${\bf n}$ is the external unit normal vector to $\ptl \Omega$,
\[
{\bf v} = (ik v, \ptl_y v),
\qquad
{\bf w} = (-ik w, \ptl_y w)
\]
are the vector--functions (the pair contains $x$- and $y$-components), ``$\cdot$'' denotes the scalar product.

The proof of (\ref{eqA103}) can be obtained in a straightforward manner by applying Gauss--Ostrogradsky theorem.

If fields $v$ and $w$ are defined on the whole surface $\bs$ and are generated by a finite set of point sources,
then (\ref{eqA103}) can be applied to surface $\bs$ and it reads as
\begin{equation}
\int \limits_\bs [f w - h v] ds = 0.
\labell{eqA103a}
\end{equation}

{\bf Proposition 1.} Directivity $V_l (\theta)$ can be represented as follows:
\begin{equation}
V_l (\theta) = 1 + 2 \sum_{n = l+1}^{\infty} e^{
i k (x_n-x_l) \theta^2/2 }
\int \limits_{-\infty}^0 w_l (x_n + 0, y, 1_*) e^{-i k \theta y} dy.
\labell{eq0506}
\end{equation}

\begin{figure}[ht]
\centerline{\epsfig{file=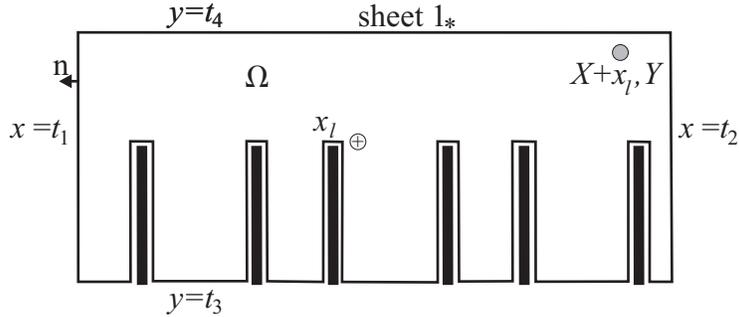}}
\caption{Domain $\Omega$}
\label{fig06a}
\end{figure}

{\bf Proof.}
Consider domain $\Omega$ shown in Fig.~\ref{fig06a}.
Substitute functions $v(x,y) = w_l(x,y,1_*)$, $w(x,y) = -G(x_l+X- x, y-Y)$ for some $X,Y$ , $Y > 0$
into
(\ref{eqA103}). Note that
\[
L' [w] = - \delta(x-(X+ x_l))\delta(y-Y),
\qquad
L [v] = \delta(x- (x_l+0))\delta(y),
\]
thus the right-hand side of (\ref{eqA103}) is equal to $ w_l (x_l + X , Y, 1_*)-G(X, Y)$.

Consider the limits $t_1, t_3 \to -\infty$, $t_2 , t_4 \to \infty$. Due to the limiting absorption
principle,
the integrals along the lines  $x = t_1, t_2$ and $y = t_3, t_4$ vanish, and in the left-hand side
of (\ref{eqA103}) one can leave only the integrals along the lines $x = x_n \pm 0$, i.\ e.\  the integrals over the shores of the cuts.

Take the limit $X \to \infty$ with $Y = X \theta$. The Green's function $G$ can be approximated as
\begin{equation}
G(X+x_l - x_n, y -X \theta) \approx \sqrt{\frac{k}{2 \pi i X}}
\exp \left\{
ik \left( X\frac{\theta^2}{2} - y \theta + (x_n - x_l) \frac{\theta^2}{2}
\right)
\right\}.
\labell{eqA104}
\end{equation}
Substituting (\ref{eqA104}) into (\ref{eqA103}) and taking into account that
\[
w_l (x_n +0 , y, 1_*) = - w_l (x_n -0 , y, 1_*)
\]
(following from (\ref{eq0403a})), get (\ref{eq0506}).

\subsection{Edge values}
\label{sec06}

Here, again, $u$ is the solution of the parabolic diffraction problem on $\bs$
formulated in Section~\ref{sec03} (with the incident plane  wave on sheet~1).
Introduce the edge values of  $u$ as follows:
\begin{equation}
C_{n, \nu_*} = u(x_n - 0,0,\nu_*).
\labell{eq0404}
\end{equation}
The notation $x_n - 0$ shows that the edge values are taken at the point located to the {\em left\/} of $x_n$, i.\ e.\
\[
C_{n, \nu_*} = \lim_{\epsilon \to 0} u(x_n - \epsilon,0,\nu_*), \qquad \epsilon>0.
\]
The points at which the edge values are taken are shown in Fig.~\ref{fig05}.

\begin{figure}[ht]
\centerline{\epsfig{file=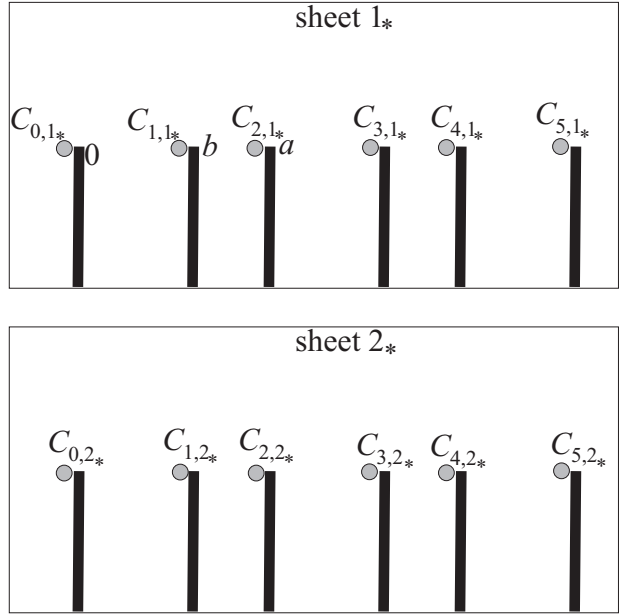}}
\caption{Positions of the points at which the edge values are taken}
\labell{fig05}
\end{figure}

{\bf Proposition 2.} The edge values $C_{n, \nu_*}$
and the directivities of the edge Green's functions $V_n$ are linked by the relation
\begin{equation}
C_{n,1_*} - C_{n,2_*} = e^{ -ikx_n \theta_{\rm in}^2/2 }
V_{1-n} (\theta_{\rm in}).
\labell{eq0701}
\end{equation}

{\bf Proof.}
Formula (\ref{eq0701}) is the parabolic analogue of the reciprocity principle. Instead of
a plane incident wave one can consider the field generated by a point source located at the point
$(-X, \theta_{\rm in} X)$ as $X\to +\infty$. The amplitude of the source should compensate the decay of
$G$ with $X$. I.\ e.
\begin{equation}
u(x,y,\nu_*) = \lim_{X \to \infty} u_X (x,y,\nu_*),
\labell{eq0702}
\end{equation}
where $u_X (x,y,\nu_*)$ obeys the following equation on $\bs$
\begin{equation}
L[u_X](x,y,\nu_*) =  A(x) \delta(x+X) \delta(y - X \theta_{\rm in})\delta_{\nu,1},
\labell{eq0703a}
\end{equation}
\[
A(X) = \sqrt{\frac{2 \pi i X}{k}} e^{ -ikX \theta_{\rm in}^2/2},
\]
it decays as $|y|\to \infty$
and is bounded near the branch points.

Consider the function
\begin{equation}
w (x,y,\nu_*) = u_X (b-x,y,\nu_*).
\labell{eq0704a}
\end{equation}
This function has the same branch points as $u$, i.\ e. it is single-valued on $\bs$.
It obeys equation
\begin{equation}
L'[w](x,y,\nu_*) = A(X)
\delta(x-(X+b)) \delta(y - X \theta_{\rm in}) \delta_{\nu,1},
\labell{eq0703b}
\end{equation}
moreover
\begin{equation}
C_{n,1_*} - C_{n,2_*} = w(x_{1-n}+0 , 0, 1_*) - w(x_{1-n}+0 , 0, 2_*).
\labell{eq0704b}
\end{equation}

Now take $v = w_{1-n} (x,y,\nu_*)$ and apply (\ref{eqA103a}) to $v$ and $w$. Consider the limit
$X \to \infty$ and take into account that $b-x_{1-n}=x_n$. By definition (\ref{eq0501}),
\begin{equation}
\lim_{X\to \infty} w_{1-n} (X+b , X\theta_{\rm in},1_*)
\sqrt{\frac{2\pi i(X+x_n)}{k}} \exp\left\{ -\frac{ik}{2} \frac{(X\theta_{\rm in})^2}{X+x_n} \right\}
= V_{1-n} (\theta_{\rm in}).
\labell{eq0705}
\end{equation}
Approximating the exponential factor as
\[
\exp\left\{ -\frac{ik}{2} \frac{(X\theta_{\rm in})^2}{X+x_n} \right\}
\approx
e^{ -ikX \theta_{\rm in}^2/2 }
e^{  ik x_n \theta_{\rm in}^2/2},
\]
obtain (\ref{eq0701}) $\square$.


\section{Embedding formulae for the field and for the scattering coefficients}
\subsection{The embedding formula for the field}
\label{sec07}

{\bf Proposition~3}. The $y$-derivative of $u$ obeys equation
\begin{equation}
L\left[ \ptl_y u \right] (x,y,\nu_*) = (-1)^{\nu-1}
\sum_{n = -\infty}^{\infty} (C_{n,1_*} - C_{n,2_*}) \delta(x-(x_n + 0)) \delta(y).
\labell{eq0405}
\end{equation}

{\bf Proof.} Function $\ptl_y u$ obeys equation (\ref{eq0302})
everywhere on $\bs$ except maybe the branch points
since operator $\ptl_y$ commutes with $L$. To study the singularities at the branch points apply
integration by parts to (\ref{eq0313}):
\begin{equation}
\ptl_y u(x, y, \nu_*) =  \int \limits_{0}^{\infty} \ptl_{y'} u(x_n,y', \nu_*) G(x-x_n , y-y') dy' +
\labell{eq0406}
\end{equation}
\[
\int \limits_{-\infty}^0 \ptl_{y'} u(x_n-0,y', (3-\nu)_*) G(x-x_n , y-y') dy' +
\]
\[
(C_{n,\nu_*} - C_{n,(3-\nu)_*}) G(x-x_n , y).
\]
Since $G$ is the field produced by a unit point source, conclude that the point $(x_n + 0, 0, \nu_*)$
contains a point source with the amplitude $C_{n,\nu_*} - C_{n,(3-\nu)_*}$, which
coincides with (\ref{eq0405}) $\square$.

Now we are ready to derive the embedding formula for the field.
Apply operator
\begin{equation}
H = \ptl_y + ik\theta_{\rm in}
\labell{eq0601}
\end{equation}
to the wave field $u(x,y,\nu_*)$.
Function $H[u](x,y)$ obeys equation
\begin{equation}
L[H[u]](x,y,\nu_*) = (-1)^{\nu-1}
\sum_{n = -\infty}^{\infty} (C_{n,1_*} - C_{n,2_*}) \delta(x-(x_n + 0)) \delta(y)
\labell{eq0602}
\end{equation}
on $\bs$ and contains no incident wave since $H[u_{\rm in}]=0$.
Thus, we can conclude that
\begin{equation}
H[u](x,y,\nu_*) = H[u_{\rm sc}](x,y,\nu_*) =
\sum_{n = -\infty}^{\infty} (C_{n,1_*} - C_{n,2_*}) w_n (x,y, \nu_*).
\labell{eq0603}
\end{equation}
The last identity is based implicitly on uniqueness of the field on $\bs$
(both sides of (\ref{eq0603}) contain no wave components coming from large $|y|$, and they
have the same set of sources, thus it follows from uniqueness that the fields
on the right and on the left are equal).
The uniqueness for the parabolic equation on a branched surface is
established in \cite{ShaninMatOsnovy} for imaginary $k$. Note that for real $k$
the uniqueness theoretically can brake due to trapped modes. This issue is not discussed in the
current paper.

Combining (\ref{eq0603}) with (\ref{eq0701}), obtain
\begin{equation}
 H[u_{\rm sc}](x,y,\nu_*) =
\sum_{n = -\infty}^{\infty} e^{ -ikx_n \theta_{\rm in}^2/2  }
V_{1-n} (\theta_{\rm in}) w_n (x,y, \nu_*),
\labell{eq0604}
\end{equation}
which is the embedding formula for the fields.

Next, the embedding formula for the fields will be converted into
the representations of the coefficients $R^{\nu}_n$, $T^{\nu}_n$
in terms of the directivities $V_n (\theta)$, i.~e.\  into embedding
formulae for the scattering coefficients. After that the coefficients
$V_n (\theta)$ will be estimated and, thus, the initial problem will be solved.


\subsection{Representation for $R_n^{\nu}$}
\label{sec08}

{\bf Proposition 4.} The values $R_n^{\nu}$ can be represented as follows:
\begin{equation}
R^1_n =
\frac{
V_0(\psi_n) V_1 (\theta_{\rm in}) + V_1(\psi_n) V_0(\theta_{\rm in})
e^{ ikb(\psi_n^2 - \theta_{\rm in}^2) /2 }
}{
ik a \psi_n (\psi_n + \theta_{\rm in}) ,
}
\labell{eq0801}
\end{equation}
\begin{equation}
R^2_n = - R^1_n.
\labell{eq0802}
\end{equation}

\begin{figure}[ht]
\centerline{\epsfig{file=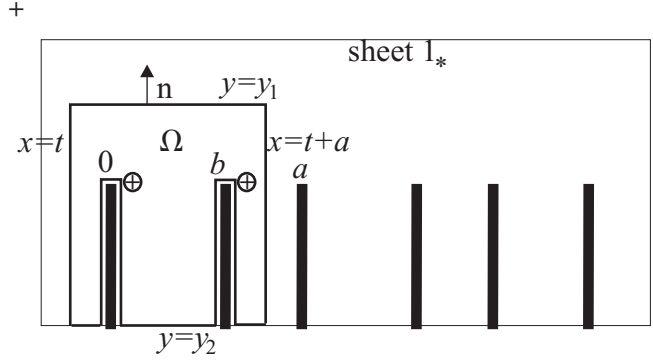}}
\caption{Domain for derivation of (\ref{eq0801})}
\label{fig06}
\end{figure}

{\bf Proof.}
Consider domain $\Omega$ shown in Fig.~\ref{fig06}. Parameters $t$ and $y_1$ are fixed,
and the limit $y_2 \to -\infty$ is taken.
Apply Green's theorem (\ref{eqA103}) in $\Omega$.
Take $v = H[u](x,y,1_*)$ in the form (\ref{eq0604}), and
\begin{equation}
w = e^{ ik x \psi^2_l/2 - i k y \psi_l }.
\labell{eq0803}
\end{equation}

Due to periodicity, the integrals over the lines $x = t$ and $x = t+a$ in the left-hand side of (\ref{eqA103})
compensate each other. Due to the limiting absorption principle, the integral over the segments
with $y = y_2$ tend to zero as $y_2 \to -\infty$.

Consider the integral in the left-hand side of (\ref{eqA103}) over the segment $y = y_1$, $t< x<t+a $.
Function $u_{\rm sc}$ can be represented as (\ref{eq0305a}).
The terms of (\ref{eq0305a}) provide Fourier decomposition on the segment. Thus, due to orthogonality of the
exponentials, only the term with $n = l$ results in a non-zero integral. Finally, (\ref{eqA103}) reads as
\begin{equation}
i k a \psi_l (\psi_l + \theta_{\rm in})R^1_l- 2 \sum_{n=-\infty}^{\infty}
e^{ -i k x_n (\theta_{\rm in})^2/2 } V_{1-n}(\theta_{\rm in})
\times
\labell{eq0804}
\end{equation}
\[
\int \limits_{-\infty}^0 \left(
w_n(+0 , y, 1_*) + e^{i k b \psi_l^2/2} w_n(b+0 , y , 1_*)
\right) e^{- i k \psi_l y} dy =
\]
\[
V_1 (\theta_{\rm in}) + e^{
i k b (\psi_l^2 - (\theta_{\rm in})^2) /2
}
V_0 (\theta_{\rm in}) .
\]

This identity can be transformed into (\ref{eq0801}) by taking into account (\ref{eq0506}),
(\ref{eq0504}), (\ref{eq0306}).

Relation (\ref{eq0802}) follows from (\ref{eq0310}) $\square$.

\subsection{Representation for $T_n^{\nu}$, $n \ne 0$}
\label{sec09}

{\bf Proposition 5.} Transmission coefficients $T^\nu_n$ for $n\ne 0$ can be represented as follows:
\begin{equation}
T^1_n =
\frac{
-V_0(\psi_n) V_1 (\theta_{\rm in}) + V_1(\psi_n) V_0(\theta_{\rm in})
e^{ ikb(\psi_n^2 - \theta_{\rm in}^2) /2 }
}{
ik a \psi_n ( \theta_{\rm in}-\psi_n ) ,
}
\labell{eq0901}
\end{equation}
\begin{equation}
T^2_n = - T^1_n.
\labell{eq0902}
\end{equation}

{\bf Proof.}
The proof of (\ref{eq0901}) is similar to the proof of (\ref{eq0801}).
Let be $u' = H[u]$ on $\bs$. Make a reflection
\begin{equation}
v(x,y,\nu) = u'(x,-y, \nu)
\labell{eq0903}
\end{equation}
and substitute it into (\ref{eqA103}) as $v$. Take (\ref{eq0803}) as $w$. Note that due to the symmetry
(\ref{eq0403}), (\ref{eq0604}) reads as
\begin{equation}
u' (x,y,\nu_*) =
\sum_{n = -\infty}^{\infty} e^{ -ikx_n \theta_{\rm in}^2/2  }
(-1)^{n+1} V_{1-n} (\theta_{\rm in}) w_n (x,y, \nu_*).
\labell{eq0904}
\end{equation}
Thus, (\ref{eq0804}) reads as
\begin{equation}
i k a \psi_l (-\psi_l + \theta_{\rm in})T^1_l- 2 \sum_{n=-\infty}^{\infty} (-1)^{n+1}
e^{ -i k x_n (\theta_{\rm in})^2/2 } V_{1-n}(\theta_{\rm in})
\times
\labell{eq0905}
\end{equation}
\[
\int \limits_{-\infty}^0 \left(
w_n(+0 , y, 1_*) + e^{i k b \psi_l^2/2} w_n(b+0 , y , 1_*)
\right) e^{- i k \psi_l y} dy =
\]
\[
-V_1 (\theta_{\rm in}) + e^{
i k b (\psi_l^2 - (\theta_{\rm in})^2) /2
}
V_0 (\theta_{\rm in}) .
\]

This results in (\ref{eq0901}). Relation (\ref{eq0902}) follows from (\ref{eq0310}) $\square$.

\subsection{Representation of $T^\nu_0$. Extended directivities}
\label{sec10}

It is necessary to find coefficient $T^\nu_0$ since the main
task of the paper is to study the coefficients of direct transmission and mirror
reflection $\tilde T_{m,m}$, $\tilde R_{m,m}$,
into which $T^\nu_0$ contributes via (\ref{eq0308}), (\ref{eq0309}).
The form of (\ref{eq0905}) makes it clear why the formula does not work for $n = 0$: the factor
$(- \psi_l+\theta_{\rm in})$ is equal to zero for $l=0$. Thus, a regularization of (\ref{eq0901}) is needed.
However, such a regularization is not a simple task, since physically one cannot vary $\psi_n$ without
changing $\theta_{\rm in}$ (they are linked by the Floquet principle). So, a more elaborate
approach is necessary.

{\bf Proposition 6.}
The regularization is given by the formula
\begin{equation}
T^1_0 = \lim_{\phi \to \theta_{\rm in}}
\frac{V_0 (\theta_{\rm in}, \phi) V_1 (\theta_{\rm in}) -
V_1 (\theta_{\rm in}, \phi) V_0 (\theta_{\rm in})}{i k a \theta_{\rm in}(\phi- \theta_{\rm in})}
\labell{eq1001}
\end{equation}
where $V_n (\theta, \phi)$ are the extended directivities given by
\begin{equation}
V_l(\theta , \phi) = \lambda_{l-1} +
2 \sum_{n=l+1}^{\infty}
\lambda_{n-1} e^{ ik (x_{n}-x_l)\theta^2/2 }
\int \limits_{-\infty}^0 w_l(x_n + 0 , y , 1_* ) e^{- i k\phi y} dy,
\labell{eq1002}
\end{equation}
\begin{equation}
\lambda_n = \lambda_n(\theta, \phi) =\left\{
\begin{array}{ll}
1 , & \mbox{even }n , \\
\exp\{ i k (a-b)( \phi^2- \theta^2) / 2 \}, &  \mbox{odd }n.
\end{array}
\right.
\labell{eq1002a}
\end{equation}
The coefficient $T^2_0$ is given by
\begin{equation}
T^2_0 = 1 - T^1_0.
\labell{eq1003}
\end{equation}

The proof of the  proposition in rather complicated and it is given in Appendix~A.

Obviously, the connection between the directivities introduced as (\ref{eq0506}) and
(\ref{eq1002}) is as follows:
\begin{equation}
V_l(\theta) = \lim_{\phi \to \theta} V_l (\phi , \theta).
\labell{eq1004}
\end{equation}

Introduce the derivatives of $V_n(\theta, \phi)$ with respect to the second argument:
\begin{equation}
V'_n (\theta) = \ptl_\phi V_n (\theta, \phi) |_{\phi = \theta}.
\labell{eq1005}
\end{equation}
By applying l'Hospital's rule, (\ref{eq1001}) can be written in the form
\begin{equation}
T^1_0 =  \frac{V'_0 (\theta_{\rm in}) V_1(\theta_{\rm in}) -
V'_1 (\theta_{\rm in}) V_0(\theta_{\rm in})}{i k a \theta_{\rm in}}.
\labell{eq1006}
\end{equation}


\section{Estimation of $V_0(\theta, \phi)$ and $V_1(\theta, \phi)$ for small $(a-b)/a$}
\subsection{Asymptotics of $V_0(\theta_{\rm in}, \phi)$ }
\label{sec11}

Consider the case of small slits between the hard screen and the walls of the waveguide:
\begin{equation}
\eps \equiv \frac{a-b}{a} \ll 1.
\labell{eq1101}
\end{equation}
Our aim is to estimate $V_0(\theta_{\rm in}, \phi)$ (and, in the next section,
$V_1(\theta_{\rm in}, \phi)$) in this case. The resulting expressions will be put into (\ref{eq0801}) and (\ref{eq1006}).
Then approximations $\tilde T_{m,m}$ and $\tilde R_{m,m}$ will be found from (\ref{eq0210}) and (\ref{eq0211}).

\begin{figure}[ht]
\centerline{\epsfig{file=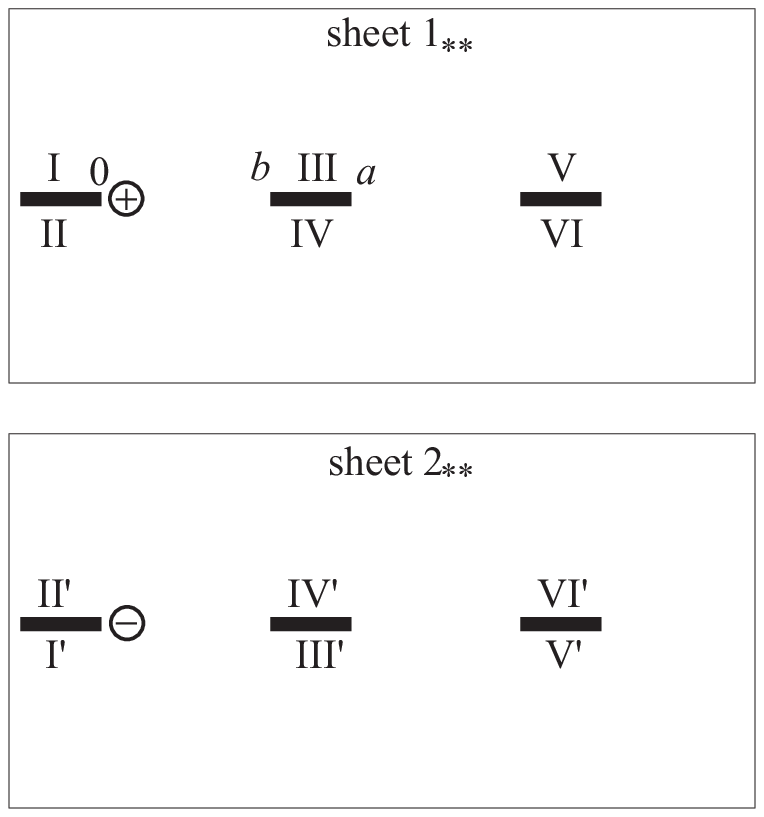}}
\caption{The third scheme of $\bs$}
\label{fig08}
\end{figure}

Cut surface $\bs$ into two sheets
differently from Fig.~\ref{fig02} and Fig.~\ref{fig03}. The scheme of the
 new cut of the surface is shown in Fig.~\ref{fig08}. New notations of the sheets are $(x,y,\nu_{**})$, $\nu = 1,2$.
The cuts are conducted now along the segments $x_{2n-1} < x < x_{2n}$. The upper half-plane of sheet~$1_{**}$
is the upper half-plane of sheet~1, i.\ e.
\begin{equation}
w_0 (x,y,\nu_{**}) = w_0 (x,y,\nu) \qquad \mbox{for } y > 0,
\labell{eq1102}
\end{equation}
\begin{equation}
w_0 (x,y,\nu_{**}) = -w_0 (x,y,\nu) \qquad \mbox{for } y < 0,
\labell{eq1103}
\end{equation}
and
\begin{equation}
w_0(x_n+0,y,\nu_*) = (-1)^n w_0(x_n,y,\nu_{**}) \qquad \mbox{for } y < 0.
\labell{eq1104}
\end{equation}

Rewrite expression (\ref{eq1002}) as follows:
\begin{equation}
V_0 (\theta, \phi) = \lambda_1 + 2 \sum_{n = 1}^{\infty} \lambda_{n-1} (-1)^n
e^{ i k x_{n} \theta^2/2 }
\int \limits_{-\infty}^0 w_0 (x_n, y, 1_{**}) e^{-i k \phi y} dy .
\labell{eq1105}
\end{equation}
Transform the sum in the right:
\begin{equation}
\lambda_1^{-1} V_0 (\theta, \phi) =
1+
2 \sum_{n = 1}^{\infty}
e^{ i k x_{2n} \theta^2/2  }
\times \qquad \qquad  \qquad  \qquad \qquad  \qquad \qquad  \qquad
\labell{eq1106}
\end{equation}
\[
\int \limits_{-\infty}^0 \left( w_0 (x_{2n} , y, 1_{**})
- e^{-ik(a-b)\phi^2/2} w_0 (x_{2n-1} , y, 1_{**})
\right) e^{-i k \phi y} dy ,
\]

Applying Green's theorem
and  taking $v = w_0(x,y)$, $w = e^{-i k \phi y + i k x \phi^2 / 2}$
rewrite (\ref{eq1106}) in the form
\begin{equation}
\lambda_1^{-1} V_0 (\theta, \phi) = 1- \frac{1}{ik}
\sum_{n = 1}^{\infty}
e^{ i k x_{2n} \theta^2/2  } \times \qquad \qquad \qquad \qquad \qquad \qquad
\labell{eq1109a}
\end{equation}
\[
\qquad \qquad \qquad
\int \limits_{x_{2n-1}}^{x_{2n}} \ptl_y w_0(x, -0, 1_{**}) e^{i k (x-x_{2n}) \phi^2/2} dx.
\]
Under the assumption $(k(a-b) \theta_{\rm in}) \ll 1$ make an approximation $\lambda_1 = 1$,  and
\begin{equation}
V_0 (\theta, \phi) \approx 1- \frac{1}{ik}
\sum_{n = 1}^{\infty}
e^{ i k x_{2n} \theta^2/2}
\int \limits_{x_{2n-1}}^{x_{2n}} \ptl_y w_0(x, -0, 1_{**})  dx
\labell{eq1109b}
\end{equation}
Here we set $e^{i k (x-x_{2n}) \phi^2/2} \approx 1$.

Consider sheet~$1_{**}$. Each cut $y < 0$, $x_{2n-1}<x<x_{2n}$
is a scatterer
on this sheet. Due to the symmetry of $w_0$, the field is equal to zero on the
cuts, so the scatterers are of Dirichlet type.
Such scatterers are studied in Appendix~B.

Introduce the coefficients
\begin{equation}
c_n = w_0(x_{2n-1}-0 , y, 1_{**}).
\labell{eq1110}
\end{equation}
According to (\ref{eqB105}),
\begin{equation}
w_0 (x,y,1_{**}) = G(x,y) + h_{\rm D} \sum_{n=1}^{\infty} c_n G(x-x_{2n} , y),
\labell{eq1115a}
\end{equation}
\[
h_{\rm D} = - 2 \sqrt{\frac{2 i (a-b)}{\pi k}}.
\]

Estimate the values $c_n$.
According to (\ref{eq1115a}),
\begin{equation}
c_n = G(n a, 0) + h_{\rm D} \sum_{l = 1}^{n-1} c_l G((n-l) a, 0)
\labell{eq1116a}
\end{equation}
or
\begin{equation}
c_n = \sqrt{\frac{k}{2 \pi i a}}
\left(
\frac{1}{\sqrt{n}}
 + h_{\rm D} \sum_{l = 1}^{\infty} \frac{c_l}{ \sqrt{n-l} }
\right).
\labell{eq1116}
\end{equation}

(\ref{eq1116}) is an infinite system of linear equations for $c_n$. Fortunately, the
system has a convolution nature (with respect to the index), so it can be solved explicitly.
Multiply (\ref{eq1116}) by $e^{i p n}$ (${\rm Im}[p] \ge 0$) and sum over $n$ from 1 to infinity. Change
the order of summation on the right:
\begin{equation}
\hat c(p) = \sqrt{\frac{k}{2 \pi i a}} {\rm Li}_{1/2} (e^{ip}) ( 1 + h \hat c(p) ),
\labell{eq1117}
\end{equation}
where
\begin{equation}
\hat c(p) = \sum_{n = 1}^{\infty} c_n e^{i p n},
\labell{eq1118}
\end{equation}
and
\begin{equation}
{\rm Li}_{\zeta} (z) = \sum_{n=1}^{\infty} \frac{z^n}{n^{\zeta}}
\labell{eq1119}
\end{equation}
is a polylogarithm function \cite{AbramowitzStegun}. Solving (\ref{eq1118}) obtain
\begin{equation}
\hat c(p) = \sqrt{\frac{k}{2 \pi i a}} \frac{{\rm Li}_{1/2}(e^{ip})}{
1 + 2 \sqrt{\eps} {\rm Li}_{1/2}(e^{ip}) / \pi
}.
\labell{eq1120}
\end{equation}
Taking into account (\ref{eq1109b}) and (\ref{eqB111}), obtain
\begin{equation}
V_0 (\theta, \phi) \approx 1 + h_{\rm D} \hat c (ka \theta^ 2/ 2) =
\left(
1 + 2 \frac{\sqrt{\eps}}{\pi} {\rm Li}_{1/2} (e^{ i k a \theta^2 /2 } )
\right)^{-1} .
\labell{eq1121}
\end{equation}

Obviously,
\[
 V'_0 (\theta) = 0.
\]


\subsection{Estimation of $V_1 (\theta, \phi)$}
\label{sec12}

Again, consider the field $w_1$ on surface $\bs$ cut according to Fig.~\ref{fig08}.
Note that due to the symmetry of the sources for $w_1$ the segments $y = 0$, $x_{2n-1}<x<x_{2n}$
can be considered as Neumann scatterers on sheet~$1_{**}$. Thus, the description developed in
Appendix~B for small Neumann segments is valid for this sheet.

Expression for $V_1 (\theta, \phi)$ has form
\begin{equation}
V_1 (\theta, \phi) = 1 + 2 \sum_{n = 2}^{\infty} \lambda_{n-1} (-1)^n
e^{ i k ( x_{n} -x_1 )  \theta^2/2 }
\int \limits_{-\infty}^0 w_1 (x_n, y, 1_{**}) e^{-i k \phi y} dy .
\labell{eq1201}
\end{equation}
Using Green's theorem rewrite (\ref{eq1201}) in the form
\begin{equation}
V_1 (\theta, \phi) = 1 + 2 e^{ i k ( x_{n} -x_1 )  \phi^2/2 }
\int \limits_{-\infty}^0 w_1 (x_2, y, 1_{**}) e^{-i k \phi y} dy -
\labell{eq1202c}
\end{equation}
\[
\phi \lambda_1 \sum_{n=2}^{\infty}
e^{ i k x_{2n} \theta^2/2  }
\int \limits_{x_{2n-1}}^{x_{2n}}  w_1(x, -0, 1_{**}) e^{i k (x-x_{2n}) \phi^2/2} dx
\]

Find function $w_1(x,y,1_{**})$ in the strip $x_2< x <x_3$. Obviously,
\begin{equation}
w_1 (a,y,1_{**}) = \mbox{sign}(y) G(a-b,y).
\labell{eq1202a}
\end{equation}
Then, for $a<x<a+b$
\[
w_1 (x,y,1_{**}) = \int \limits_{-\infty}^{\infty}
w_1 (a,y',1_{**}) G(x-a , y-y') dy'
\]
Repeat the argument from Appendix~B and expand $G(x-a , y - y')$ as a Taylor series in $y'$:
\[
w_1 (x,y,1_{**}) \approx \int \limits_{-\infty}^{\infty}
w_1 (a,y',1_{**}) \left( (G(x-a , y)
- y' \ptl_y G(x-a , y) \right)
dy'.
\]
Note that the first term in the parentheses (the monopole term) yields zero. Thus,
\begin{equation}
w_1 (x,y,1_{**}) \approx f_0 \ptl_y G(x-a , y) ,
\labell{eq1203c}
\end{equation}
\[
f_0 = - 2 \int \limits_0^{\infty} y G(a-b, y) dy = -\sqrt{\frac{2 i (a-b)}{\pi k}}
\]
The field (\ref{eq1203}) near the point $(x_3, 0)$ has ``dipole'' nature (see Appendix~B),
i.\ e.\ it can be approximated as $c y$ on a large segment $x = x_3$, $-\sqrt{a/(k\eps)}<y<\sqrt{a/(k\eps)}$.

Redefine the values
\begin{equation}
c_n = \ptl_y w_1 (x_{2 n +1 }-0, 0, 1_{**}).
\labell{eq1204}
\end{equation}
According to the approximations developed in Appendix~B, for $y > a$ one can write
\begin{equation}
w_1 (x , y, 1_{**}) \approx
f_0 \ptl_y G(x-a , y) + h_{\rm N} \sum_{n = 1}^\infty
c_n \ptl_y G(x-a(n+1), y),
\labell{eq1205}
\end{equation}
\[
h_{\rm N} = -\frac{2}{3} \
\left( \frac{a-b}{k} \right)^{3/2} \sqrt{\frac{2}{\pi i}}
\]
A system of equations for $c_n$ has form
\begin{equation}
c_n = f_1 \left( \frac{f_0}{n^{3/2}} + h_{\rm N} \sum_{l=1}^{n-1} \frac{c_l}{(n-m)^{3/2}} \right),
\labell{eq1206}
\end{equation}
\[
f_1 = n^{3/2} \ptl^2_y G(a n , 0) = \left(\frac{k}{a}\right)^{3/2} \sqrt{\frac{i}{2 \pi}}.
\]
Introduce the function $\hat c(p)$ by formula (\ref{eq1118}). Solve (\ref{eq1206}) by Fourier transform:
\begin{equation}
\hat c(p) = \frac{f_0 f_1 {\rm Li}_{3/2} (e^{ip})}{1 - f_1 h_{\rm N} {\rm Li}_{3/2} (e^{ip})}.
\labell{eq1207}
\end{equation}

According to (\ref{eq1202}), (\ref{eq1202a}), and (\ref{eqB112}),
\begin{equation}
V_1 (\theta, \phi) = 1 - 2 e^{ i k ( a-b )  \phi^2/2 }
\int \limits_{-\infty}^0 G (a-b, y) e^{-i k \phi y} dy +
\labell{eq1202}
\end{equation}
\[
ik \phi \lambda_1 h_{\rm N}
e^{ i k a \theta^2/2  } \hat c(k a \theta^2/2)
\]
An elementary estimation of the first two terms yields
\[
1 - 2 e^{ i k ( a-b )  \phi^2/2 }
\int \limits_{-\infty}^0 G (a-b, y) e^{-i k \phi y} dy =
\sqrt{\frac{2 (a-b) k}{\pi i}} \phi + O(k a \theta^2).
\]
Thus, this value is of order $\phi \sqrt{ka} \sqrt{\eps}$.
Make an estimation of the third term of (\ref{eq1202}):
\[
ik \phi \lambda_1 h_{\rm N}
e^{ i k a \theta^2/2  } \hat c(k a \theta^2/2) \sim \phi \sqrt{ka} \eps^2.
\]
One can see that the third term contains an additional small factor $\eps^{3/2}$ and can be neglected.
Thus, a usable approximation of $V_1 (\theta, \phi)$ is
\begin{equation}
V_1 (\theta, \phi) \approx \sqrt{\frac{2 (a-b) k}{\pi i}} \phi,
\labell{eq1203}
\end{equation}
and
\begin{equation}
V'_1 (\theta) \approx \sqrt{\frac{2 (a-b) k}{\pi i}}.
\labell{eq1203a}
\end{equation}


\section{Asymptotical and numerical results}
\label{sec13}

\begin{figure}[ht]
\centerline{\epsfig{file=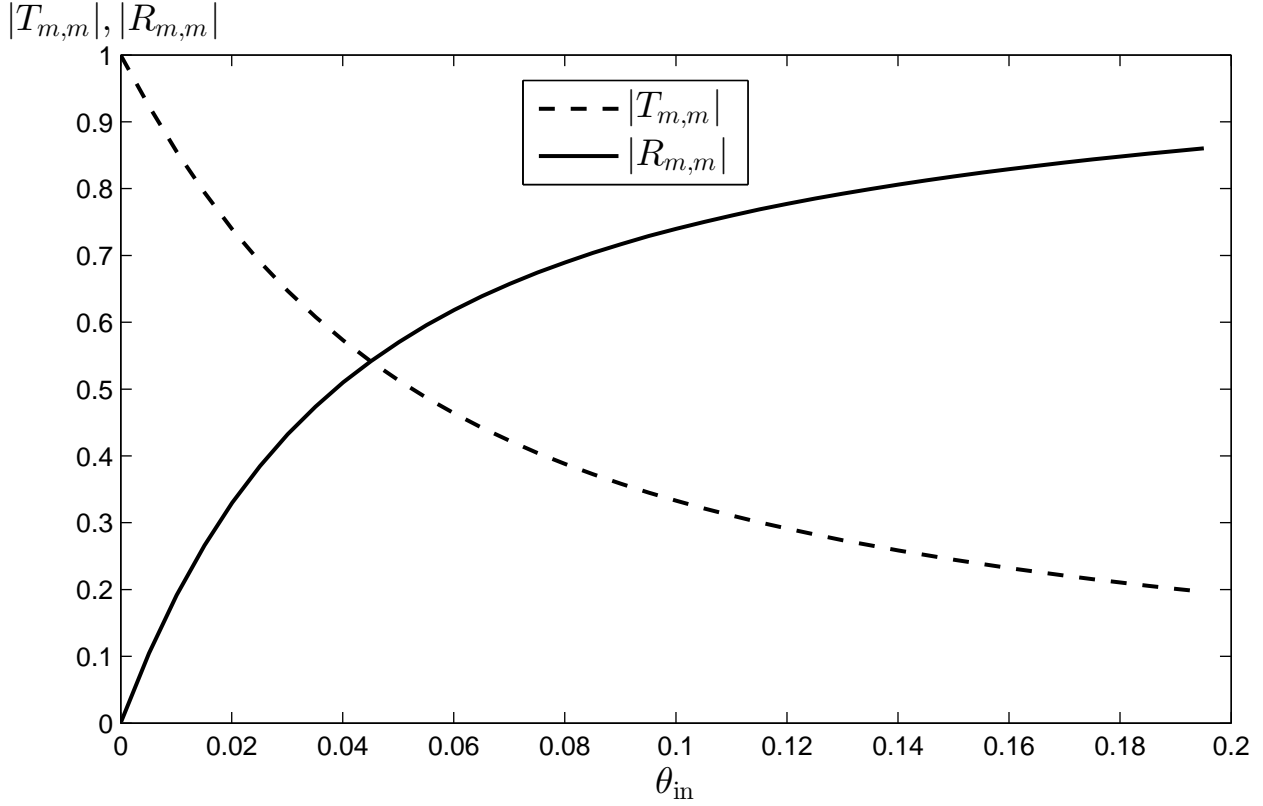}}
\caption{Transmission and reflection coefficients $\tilde T_{m,m}$, $\tilde R_{m,m}$}
\label{fig09}
\end{figure}

First, estimate the straight transmission and mirror reflection coefficients $\tilde T_{m,m}$
and $\tilde R_{m,m}$  for small $(a-b)/a$.
The estimation can be obtained by applying formulae (\ref{eq0308}), (\ref{eq0309}),
(\ref{eq0310}), (\ref{eq0801}), (\ref{eq1006}). One should substitute the approximations
(\ref{eq1121}) and (\ref{eq1203}) into these formulae.
A graph of $|\tilde T_{m,m} (\theta_{\rm in})|$, $|\tilde R_{m,m} (\theta_{\rm in})|$
is shown in Fig.~\ref{fig09}. The values of parameters taken for these computations are as follows:
$ak = 100$, $(a-b)/a = 0.05$. Variation of $\theta_{\rm in}$  can be performed physically by slight
variation of temporal frequency near the threshold value. The index of the incident
waveguide mode is equal to $m = 31$. One can see that $|\tilde T_{m,m}| \to 1$ and
$|\tilde R_{m,m}| \to 0$ as $\theta_{\rm in} \to 0$. This agrees with conclusions of the consideration
developed by S.~A.~Nazarov~\cite{Nazarov2016}.

One can study (\ref{eq1121}) asymptotically for small $\theta$ (such that $\sqrt{ka} \, \theta \ll 1$).
For this, use the asymptotics of ${\rm Li}_{1/2} (e^\mu)$ for small $\mu$ \cite{AbramowitzStegun}:
\begin{equation}
{\rm Li}_{1/2} (e^{\mu}) = \sqrt{\pi} (-\mu)^{-1/2} + \dots
\labell{eq1122}
\end{equation}
Thus,
\begin{equation}
V_0 (\theta , \phi) \approx 1
\qquad \mbox{ for } \theta \gg \sqrt{\frac{\eps}{ak}},
\labell{eq1123}
\end{equation}
\begin{equation}
V_0 (\theta, \phi) \approx \frac{\theta}{2} \sqrt{\frac{\pi ak}{2 i \eps}}
\qquad \mbox{ for } \theta \ll \sqrt{\frac{\eps}{ak}},
\labell{eq1124}
\end{equation}
Using these approximations one can  obtain for small $\eps$
\begin{equation}
\tilde T_{m,m} \approx 1, \qquad |\tilde R_{m,m}| << 1 \qquad \mbox{for }\theta \ll \sqrt{\frac{\eps}{ak}},
\labell{eq1301}
\end{equation}
\begin{equation}
|\tilde T_{m,m}| \ll 1, \qquad \tilde R_{m,m} \approx 1 \qquad \mbox{for }\theta \gg \sqrt{\frac{\eps}{ak}}.
\labell{eq1302}
\end{equation}
This is the main asymptotical result of the paper. This result establishes the domain of the parameters in which
the Nazarov's anomalous transmission can happen when the gap $a-b$ is small.

\begin{figure}[ht]
\centerline{\epsfig{file=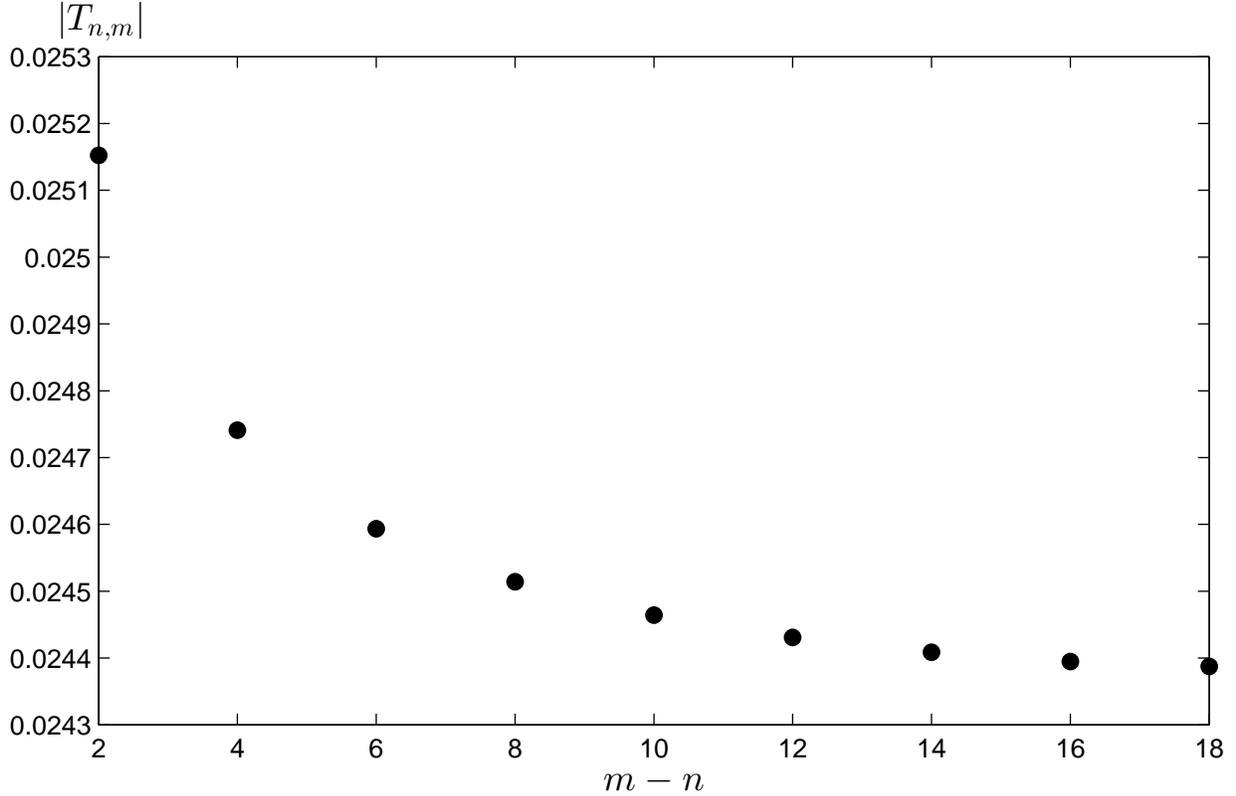}}
\caption{Transmission and reflection coefficients $\tilde T_{n,m}$, $\tilde R_{n,m}$
for the intermediate regime}
\label{fig10}
\end{figure}

The formulae obtained in the paper enable one to describe the whole process of reflection / transmission in a waveguide with a screen. Consider the case $\theta \sim \sqrt{\eps/(ak)}$, i.\ e.\ neither asymptotics (\ref{eq1301}) nor (\ref{eq1302}) is valid. Namely, take the paprameters $\theta_{\rm in} = 0.045$, $ka = 100$, $(a-b)/a = 0.05$.  Note that $\tilde T_{n,m}=-\tilde R_{n,m}$ if $m\neq n$. The values
$|\tilde T_{n,m}|$ are shown in Fig.~\ref{fig10}.

\section{Summary}

The problem of scattering by a thin Neumann wall placed in a 2D planar waveguide is studied.
The wall is irradiated by a single waveguide mode with frequency close to its cut-off.
The main aim of the paper
is to study the anomalous transmission, i.\ e.\ the phenomenon of an almost complete transmission as the
frequency of the incident mode tends to the mode's cut-off. The problems seems non-trivial when the
frequency is close to the cut-off, and the gap between the screen and the wall of the waveguide is small.

By the method of reflections the problem is first transformed into a problem of scattering by an infinite diffraction grating
composed of thin Neumann screens, and then into a problem of propagation on a branched surface having
two sheets and a periodic set of branch points of order~2.

A parabolic approximation of diffraction theory is used to describe the wave process on the
branched surface. This simplifies the consideration a lot.

The study is performed in two major steps. First, the edge Green's functions are introduced and the embedding formulae for the coefficients of reflection and transmission are derived. These formulae ((\ref{eq0801}),
(\ref{eq0901}), (\ref{eq1006})) express the coefficients of reflection and transmission in terms of the directivities
of the edge Green's functions.
Note that formulae (\ref{eq0801}) and (\ref{eq0901}) are rather standard (they are obtained by applying an embedding operator), while (\ref{eq1006}) is less trivial. This formula is a regularization of (\ref{eq0901}). The regularization is based on an extended formulation of the diffraction problem. This formulation is described in Appendix~A.

The second step is an estimation of the edge Green's functions $V_0 (\theta)$, $V_1(\theta)$ and their
extended versions $V_0(\theta, \phi)$, $V_1 (\theta, \phi)$. Such an estimation is done for a small gap, i.~e.\
for small $(a-b)/a$. The estimation is based on the description of scattering by small Dirichlet and Neumann
segments developed in Appendix~B. It is shown that a Dirichlet segment under certain conditions can be considered
as a point (monopole) scatterer having known strength. Similarly, a small Neumann segment can be considered as a dipole scatterer of a known strength. The estimation of $V_0$ and $V_1$ is then built from a description of an infinite periodic array of such scatterers.

Finally, the estimations of $V_0$ and $V_1$ are substituted into the embedding formulae, and the coefficients
of transmission and reflection for the initial waveguide problem are calculated. The main result related to the
anomalous transmission is (\ref{eq1301}), (\ref{eq1302}).

Note that the same consideration can be applied to the problem of diffraction by a thin Dirichlet screen. Also,
a small screen can be studied instead of a small gap.
In this case it can be shown that near the cut-off frequency the reflection coefficient $R_{m,m}$ is close to -1.


\section*{Acknowledgements}

The authors are
grateful to participants of the seminar on wave diffraction held in S.Pb. branch of Steklov
Mathematical Institute of RAS (the chairman is Prof. V.\ M.\ Babich) for interesting discussions.

The work is supported by the grants RFBR 14-02-00573 and Scientific Schools-283.2014.2.


\section{Appendix A. A proof of (\ref{eq1001})}

To prove (\ref{eq1001}) introduce a new formulation of a diffraction problem on $\bs$.
The ``old'' formulation depends on a single parameter $\theta_{\rm in}$, while the ``new''
formulation depends on two variables, $\theta_{\rm in}$ and $\phi$, i.\ e.\ it is more flexible.
The old formulation can be obtained
from the new one as a result of the limiting procedure $\phi \to \theta_{\rm in}$.
Such a limiting procedure provides a regularization for (\ref{eq0901}).

Consider the ``old'' diffraction problem. Redefine the scattered field (introduce
the new scattered field $u_{\rm scn}$):
\begin{equation}
u_{\rm scn} (x,y,1_*) = u(x,y,1_*) - u_{\rm in}(x,y),
\qquad
u_{\rm scn} (x,y,2_*) = u(x,y,2_*).
\labell{eqB01}
\end{equation}
The difference between (\ref{eqB01}) and (\ref{eq0304a}) is that
\[
u_{\rm scn}(x,y,1_*) = u_{\rm sc}(x,y,1_*) - u_{\rm in}(x,y)
\qquad
\mbox{for} \quad
y < 0.
\]
A direct consequence from (\ref{eqB01}) is that
\begin{equation}
u (x,0,1) = u_{\rm scn}(x,0,1) + e^{
-i k x \theta_{\rm in}^2/2
}
\labell{eqB01a}
\end{equation}
for $ b< x<a$.

Our aim here is finding the coefficient of the Fourier series:
\begin{equation}
T^1_0 = \frac{1}{a} \int \limits_0^a
u(x,-0,1)
e^{
i k x \theta_{\rm in}^2/2
} dx =
\labell{eqB01b}
\end{equation}
\[
\frac{1}{a} \int \limits_0^a
u_{\rm scn}(x,-0,1)
e^{
i k x \theta_{\rm in}^2/2
} dx
+ \frac{a-b}{a} .
\]

Since the total field is continuous on $\bs$, the new scattered field obeys matching conditions
\begin{equation}
u_{\rm scn} (x_n +0 , y, 2_*) - u_{\rm scn} (x_n - 0 , y, 1_*) =
 e^{ - i k x_n \theta_{\rm in}^2/2 }
e^{ - i k \theta_{\rm in} y },
\labell{eqB102}
\end{equation}
\begin{equation}
u_{\rm scn} (x_n +0 , y, 1_*) - u_{\rm scn} (x_n - 0 , y, 2_*) =-
 e^{ - i k x_n \theta_{\rm in}^2/2 }
e^{ - i k \theta_{\rm in} y }
\labell{eqB103}
\end{equation}
for $y < 0$.
These matching conditions (taken with the parabolic equation (\ref{eq0302}) for $u_{\rm scn}$)
can be put in the basis of the new formulation of a diffraction problem on
$\bs$. It is necessary also to take into account
the limiting absorption principle, i.~e.\ the field $u_{\rm scn}$ should decay as $|y|\to \infty$.
The field  should be bounded near the branch points as well.
As it is shown in \cite{ShaninMatOsnovy},
the field $u_{\rm scn}$ obeying these conditions is defined uniquely.

Now introduce an extended problem formulation. Find the field $u_{\rm ext}(x,y,\nu)$ on $\bs$ obeying parabolic
equation (\ref{eq0302}) everywhere except the cuts $x = x_n$, $n \in \mathbb{Z}$, $y < 0$ and
matching conditions on the cuts similar to (\ref{eqB102}), (\ref{eqB103}):
\begin{equation}
u_{\rm ext} (x_n +0 , y, 2_*) - u_{\rm ext} (x_n - 0 , y, 1_*) =
\lambda_n e^{ - i k x_n \theta_{\rm in}^2/2}
e^{ - i k \phi y },
\labell{eqB104}
\end{equation}
\begin{equation}
u_{\rm ext} (x_n +0 , y, 1_*) - u_{\rm ext} (x_n - 0 , y, 2_*) =-
\lambda_n e^{ - i k x_n \theta_{\rm in}^2/2 }
e^{ - i k \phi y },
\labell{eqB105}
\end{equation}
decaying as $|y| \to \infty$, and bounded near the points $(x_n, 0)$.
Parameter $\lambda_n = \lambda_n(\theta_{\rm in}, \phi)$ is defined by (\ref{eq1002a}).

The solution of the problem for $u_{\rm ext}$ is also unique.
This can be easily proven by writing down an explicit solution
for the case of $k$ having a positive imaginary part (see \cite{ShaninMatOsnovy}).
Function $u_{\rm ext}$
should obey the Floquet property
\[
u_{\rm ext} (x + a, y, \nu_*) =
e^{
- i k a \theta_{\rm in }^2/2 }
u_{\rm ext} (x , y, \nu_*).
\]

Under the condition
${\rm Im}[k] > 0$ the solution $u_{\rm ext}$ is smooth with respect to parameter
$\phi$, thus
\begin{equation}
u_{\rm scn} = \lim_{\phi \to \theta_{\rm in}} u_{\rm ext} ,
\labell{eqB106}
\end{equation}
and
\begin{equation}
T^1_0 = \frac{1}{a} \lim_{\phi \to \theta_{\rm in}}
 \int \limits_0^a
u_{\rm ext}(x,-0,1)
e^{
i k x \theta_{\rm in}^2/2
} dx
+ \frac{a-b}{a} .
\labell{eqB107}
\end{equation}

In the ``extended'' formulation (\ref{eqB104}), (\ref{eqB105}) one can choose any
right-hand side depending on parameter $\phi$ such that it tends to the right-hand side of
(\ref{eqB102}), (\ref{eqB103}) as $\phi \to \theta_{\rm in}$. Our choice is motivated by the following:

\noindent
--- The right-hand side of (\ref{eqB104}), (\ref{eqB105}) obeys the same Floquet property as~$u$.

\noindent
--- Function $e^{-i k \phi y}$ is nullified by operator $H_\phi = \ptl_y + i k \phi$, which is
close (but not equal) to the embedding operator~$H$.

\noindent
--- Parameters $\lambda_n$ are introduced to simplify the computations in the subsequent estimation of
$V_0(\theta, \phi)$ and $V_1(\theta, \phi)$. The choice $\lambda_n \equiv 1$ is also possible, but leads
to estimation of additional terms. Of course, the result is the same.

Consider the function
$
 H_{\phi}[ u_{\rm ext}](x,y,\nu_*).
$
defined on $\bs$.
This function decays as $|y| \to \infty$.
Moreover, this field is continuous on the cuts, since the discontinuity $\sim e^{-i k \phi y}$
is nullified by $H_\phi$. Thus, $H_\phi[u_{\rm ext}]$ is continuous on $\bs$, obeys (\ref{eq0302}) on the cuts, and the only singularities it can have are some singularities at the branch points.

{\bf Proposition A1.}
Field $v$ obeys an inhomogeneous parabolic equation
\begin{equation}
L [H_\phi[u_{\rm ext}]](x,y,\nu_*) =
\qquad \qquad \qquad \qquad \qquad \qquad \qquad \qquad \qquad
\labell{eqB108}
\end{equation}
\[
(-1)^{\nu-1}\sum_{n = -\infty}^{\infty} \left(
c_{n,1_*} - c_{n,2_*} + \lambda_n e^{ - i k x_n \theta_{\rm in}^2/2 }
\right)
 \delta (x - (x_n + 0)) \delta (y),
\]

where
\begin{equation}
c_{n,\nu_*} = u_{\rm ext} (x_n - 0, 0 , \nu_*).
\labell{eqB109}
\end{equation}

{\bf Proof.}
Prove that the singularity on the sheet $1_*$ near the point $x_0$ is
$\delta(x - 0)\delta(y) (c_{0,1_*} -c_{0,2_*} + 1)$. All other branch points / sheets can be considered
similarly. Take arbitrary $y$ and $0<x<b$. Taking into account (\ref{eqB103}), obtain a relation
\begin{equation}
u_{\rm ext} (x,y,1_*) =
\int \limits_{0}^{\infty}  u_{\rm ext} (0,y',1_*) G(x,y-y') dy' +
\labell{eqB110}
\end{equation}
\[
\int \limits_{-\infty}^{0} (u_{\rm ext} (-0,y',2_*) - e^{-i k \phi y})  G(x,y-y') dy'
\]
Apply operator $H_{\phi}$ to (\ref{eqB110}). Performing integration by parts, obtain
\begin{equation}
H_{\phi}[u_{\rm ext}] (x,y,1_*) =
\int \limits_{0}^{\infty} H_{\phi}[ u_{\rm ext} ] (0,y',1_*) G(x,y-y') dy'+
\labell{eqB111}
\end{equation}
\[
\int \limits_{-\infty}^{0} H_{\phi}[u_{\rm ext}] (-0,y',2_*)   G(x,y-y') dy' +
(c_{0,1_*} -c_{0,2_*} + 1) G(x,y).
\]
The last term corresponds to the singularity of the form $\delta(x - 0)\delta(y) (c_{0,1_*} -c_{0,2_*} + 1)$
$\square$.

{\bf Proposition A2.}
The combination $c_{n,1_*} -c_{n,2_*} + \lambda_n e^{-i k x_n \theta^2/2}$
is expressed in terms of the extended directivities as
follows:
\begin{equation}
c_{n,1_*} -c_{n,2_*} + \lambda_n e^{-i k x_n \theta^2/2} =
e^{
- i k x_n \theta_{\rm in}^2/2
}
V_{1-n} (\theta_{\rm in} , \phi)
\labell{eqB112}
\end{equation}

{\bf Proof.}
Note that by formulation (\ref{eqB104}), (\ref{eqB105})  function $u_{\rm ext}$
has a symmetry
\begin{equation}
u_{\rm ext} (x, y, 1_*) = -
u_{\rm ext} (x, y, 2_*).
\labell{eqB113}
\end{equation}
Then step by step follow the proof of (\ref{eq0701}) $\square$.

Substitute (\ref{eqB112}) into  (\ref{eqB108}) and apply the uniqueness:
\begin{equation}
H_\phi[u_{\rm ext}](x,y, \nu_*) = \sum_{n=-\infty}^{\infty}
e^{
- i k x_n \theta_{\rm in}^2/2
} \,
V_{1-n} (\theta_{\rm in} , \phi) \, w_n(x,y,\nu_*).
\labell{eqB114}
\end{equation}

{\bf Proposition A3.}
For $y < 0$
\begin{equation}
H_\phi[u_{\rm ext}](x,y,1) =
\sum_{n = -\infty}^{\infty} B_n
e^{-i k \psi_n y- i k x \psi_n^2/2 } ,
\labell{eqB115}
\end{equation}
\begin{equation}
B_n = \frac{1}{ \psi_n a} \left(
e^{
ik b (\psi_n^2 - \theta_{\rm in}^2)/2
}
V_1 (\psi_n) V_0 (\theta_{\rm in}, \phi)
- V_0(\psi_n) V_1 (\theta_{\rm in}, \phi)
\right)
\labell{eqB116}
\end{equation}

The proof is similar to that of (\ref{eq0901}).

Now we can prove (\ref{eq1001}). Obviously, $u_{\rm ext} (x,y,1)$ for
$y < 0$ can be represented as a sum of a regular and discontinuous component
\begin{equation}
u_{\rm ext} (x,y,1) =
u_{\rm reg} (x, y) + u_{\rm dis} (x,y).
\labell{eqB117}
\end{equation}
The regular component can be obtained by formal inversion of the operator $H_\phi$ at (\ref{eqB115}):
\begin{equation}
u_{\rm reg} (x, y) =
\sum_{n = -\infty}^{\infty} \frac{B_n e^{-i k \psi_n y- i k x \psi_n^2/2 }}{i k (\phi - \psi_n)}
.
\labell{eqB118}
\end{equation}
The singular part should be nullified by $H_\phi$, obey equation (\ref{eq0302}), and obey
(\ref{eqB104}), (\ref{eqB105}) on the cuts. The only form for such a solution is as follows:
\begin{equation}
u_{\rm dis}(x,y) = g(x) e^{
- i k \phi y - i k x \phi^2/2
}\qquad
\mbox{for }
0 < x < a,
\labell{eqB119}
\end{equation}
where
\begin{equation}
g(x) = g_1
\qquad
\mbox{for }
0 < x < b,
\labell{eqB120}
\end{equation}
\begin{equation}
g(x) = g_2
\qquad
\mbox{for }
b < x < a,
\labell{eqB121}
\end{equation}
and function $u_{\rm dis}(x,y)$ should obey the Floquet property
\[
u_{\rm dis} (x + a, y) =
e^{
- i k a \theta_{\rm in }^2/2 }
u_{\rm ext} (x , y).
\]

The constants $g_1$ and $g_2$ should obey equations following from the matching conditions
(\ref{eqB104}), (\ref{eqB105}):
\begin{equation}
(g_1 - g_2) e^{ i k a (\theta_{\rm in}^2 - \phi^2 )/2  } = 1,
\labell{eqB122}
\end{equation}
\begin{equation}
g_1 -  g_2 e^{ i k a (\theta_{\rm in}^2 - \phi^2 )/2 }  = 1.
\labell{eqB123}
\end{equation}
One can easily find the solution of (\ref{eqB122}), (\ref{eqB123}):
\[
g_1 =0 , \qquad g_2 = -e^{i k a (\phi^2 - \theta_{\rm in}^2  )/2}
\]
Substituting $u^{\rm reg}$ and $u^{\rm dis}$ into (\ref{eqB107}) and taking the
limit $\phi \to \theta_{\rm in}$  obtain (\ref{eq1001}).


\section*{Appendix B. Scattering by small segments}

Consider two auxiliary problems: scattering of a locally plane wave by a small Dirichlet segment
and scattering of a ``dipole'' wave by a Neumann segment. In both cases the scatterer is the segment
$y=0$, $0<x<q$ (in our case $q = a-b$),
and the direction of propagation of the wave is the positive $x$-direction. The problem is considered
in the parabolic approximation, i.\ e.\ the governing equation is (\ref{eq0302}).

Note that the concept of a small segment in the parabolic consideration is different from the Helmholtz case.
In the Helmholtz formulation the segment is small if it is smaller than the wavelength: $kq \ll 1$. In the parabolic
consideration there exists an angular parameter of the problem (in our case this parameter is $\theta_{\rm in}$)
which should be compared to the angular parameter associated with the segment: $\theta_{\rm seg} = (kq)^{-1/2}$.  The segment is small if $\theta_{\rm in} \ll \theta_{\rm seg}$.
The diffraction directivity of the Dirichlet segment is approximately constant
for $|\theta| \ll \theta_{\rm seg}$.

Another important parameter for the auxiliary problem is the diffraction width $d = \sqrt{q/k}$. This is the width
(in the $y$-direction) of the strip, which affects the diffraction, i.~e.\ the diffracted field ``feels'' the incident
wave only in the strip $-Nd < y < N d$ for sufficiently large $N$ (say $N\sim 10$).

Both $\theta_{\rm seg}$ and $d$ are well known in diffraction theory. The first one is (up to a constant factor)
the width
of the first diffraction maximum.
The second one is the width of the penumbral zone.

First, consider diffraction of a locally plane wave by by a Dirichlet strip. The locally plane wave is the wave
which is approximately constant on a segment $x = 0$, $- N d<  y <Nd $.

Denote the field falling on the segment by $u_{\rm iD} (x,y)$. Let this field obey equation (\ref{eq0301}).
Let $u_{\rm iD}$ be even: $u_{\rm iD}(x,y) = u_{\rm iD}(x,-y)$.
Define
\begin{equation}
c_{\rm D} = u_{\rm iD} (0, 0).
\labell{eqB101}
\end{equation}

The total field can be found explicitly from (\ref{eq0311}) (see also \cite{Shanin1}). The Dirichlet boundary condition leads to an odd
continuation of the field through the segment, thus
\begin{equation}
u(q , y) = \int \limits_0^{\infty} u_{\rm iD} (0,y') G(q , |y|-y') dy' -
\int \limits_0^{\infty} u_{\rm iD} (0,y') G(q , |y|+y') dy'.
\labell{eqB102a}
\end{equation}
Using the elementary property of the parabolic equation, namely
\[
u_{\rm iD}(d,y) = \int \limits_{-\infty}^{\infty} u_{\rm iD}(0,y') G(q,y-y') dy',
\]
transform (\ref{eqB102a}) into
\begin{equation}
u(q , y) = u_{\rm iD} -
2 \int \limits_0^{\infty} u_{\rm iD} (0,y') G(q , |y|+y') dy'.
\labell{eqB102b}
\end{equation}
The second term can be considered as the diffracted field $u_{\rm dD}(x,y)$
on the line $x = q$. This field has the following properties:

\noindent
--- $u_{\rm dD} (q, y)$ is small comparatively to $u_{\rm dD}(q, 0)$ for $y \gg d$.

\noindent
--- In the integral $u_{\rm iD} (0,y')$ can be replaced by its value at $y'=0$, i.\ e.\
by~$c_{\rm D}$. This will not affect the integral in the area $y \sim d$.

\noindent
--- The integral of $u_{\rm dD} (q, y)$ over $y$ is equal to
\begin{equation}
\int \limits_{-\infty}^{\infty} u_{\rm dD}(q+0 , y) dy \approx c _{\rm D}\, h_{\rm D},
\labell{eqB103c}
\end{equation}
\begin{equation}
h_{\rm D} = -4 \int \!\!\!\! \int_0^{\infty} G(q,y+y') dy dy'=- 2 \sqrt{\frac{2 i q}{\pi k}}.
\labell{eqB104b}
\end{equation}
Since for $x > q$ the diffracted field can be written as
\[
u_{\rm dD}(x , y) = \int \limits_{-\infty}^{\infty}
u_{\rm dD}(q, y') G(x-q, y-y') d y',
\]
an we can conclude that
\begin{equation}
u_{\rm dD}(x,y) \approx c_{\rm D}\, h_{\rm D} G(x,y).
\labell{eqB105b}
\end{equation}
This approximation is valid when $x \gg q$ and the angle of scattering is smaller than $\theta_{\rm seg}$. Formula (\ref{eqB105b})
means that a small Dirichlet segment acts as a point scatterer with force~$h_{\rm D}$. One can think about
 the scattered field as about a solution of an inhomogeneous parabolic equation with the source $c_{\rm D}\, h_{\rm D} \delta(x) \delta(y)$.

Now consider diffraction of a ``dipole'' wave by a Neumann segment. The position of the segment is the same as above. The dipole--type wave is an odd field $u_{\rm iN}(x,y)$ such that
$u_{\rm iN}(x,y) = - u_{\rm iN}(x,-y)$ and
\begin{equation}
u_{\rm iN}(0,y) \approx c_{\rm N} \, y
\labell{eqB106c}
\end{equation}
on some segment $-N d < y < Nd$.

Instead of the odd continuation (\ref{eqB102a}) one should use an even continuation.
Similarly to the consideration of the Dirichlet segment, the diffracted field for $x = q$ is as follows
:
\begin{equation}
u_{\rm dN}(q , y) \approx \pm 2 c_{\rm N} \int \limits_0^{\infty} y' G(q,|y|+y') dy',
\labell{eqB107c}
\end{equation}
where the upper sign is taken for $y > 0$, and the lower sign is taken for $y < 0$.

Find the field $u_{\rm dN} (x, d)$ for $x \gg q$:
\begin{equation}
u_{\rm dN}(x, y) = \int \limits_{-\infty}^{\infty} u_{\rm dN}(q , y')
G(x-q , y-y') dy'
\labell{eqB108d}
\end{equation}
Represent the kernel $G(x-q, y-y')$ as Taylor series with respect to $y'$
(this representation is valid on the segment of size $\sim d$). Also substitute $x-q$ by $x$:
\begin{equation}
u_{\rm dN}(x, y) \approx \int \limits_{-\infty}^{\infty} u_{\rm dN}(q , y')
\left( G(x , y) - y' \ptl_y G(x,y)  \right)  dy'.
\labell{eqB108e}
\end{equation}
The first term is equal to zero, since $u_{\rm dN} (q , y)$ is an odd function. Thus,
\begin{equation}
u_{\rm dN}(x, y) \approx  c_{\rm N}\, h_{\rm N}\ptl_y G(x,y),
\labell{eqB109c}
\end{equation}
\begin{equation}
h_{\rm N} = 4 \int \! \! \! \! \int_0^\infty y' y G(q, y + y') dy dy' = -\frac{2}{3} \
\left( \frac{q}{k} \right)^{3/2} \sqrt{\frac{2}{\pi i}}
\labell{eqB110c}
\end{equation}
The scattered wave is generated by the source having form $c_{\rm N}\, h_{\rm N} \delta(x) \delta'(y)$, where
$\delta'(y)$ is the derivative of the delta-function.

According to Green's theorem, formula (\ref{eqB105b}) for the Dirichlet strip is equivalent to
the formula
\begin{equation}
\int \limits_0^q \ptl_y u_{\rm dD} (x, +0) dx = - \int \limits_0^q \ptl_y u_{\rm dD} (x, -0) dx
\approx \i k c_{\rm D}\, h_{\rm D}.
\labell{eqB111c}
\end{equation}
Similarly, for the Neumann strip and a dipole field
\begin{equation}
\int \limits_0^q  u_{\rm dN} (x, +0) dx = - \int \limits_0^q u_{\rm dN} (x, -0) dx
\approx  i k c_{\rm N}\, h_{\rm N}.
\labell{eqB112c}
\end{equation}


\end{document}